\documentclass[12pt]{amsart}
\usepackage{amsmath,amssymb,amsfonts,mathtools}
\usepackage{mathrsfs}
\usepackage[shortlabels]{enumitem}
\setlist{leftmargin=*}
\usepackage{pb-diagram}

\newcommand{\eps}{\epsilon}

\gdef\Smu{\operatorname{Stab}^\mu}
\gdef\Stab{\operatorname{Stab}}

\newcommand{\al}{\alpha}

 \newcommand{\CK}{\mathcal K}
 \newcommand{\CH}{\mathcal H}
 \newcommand{\RR}{\mathbb{R}}

 \newcommand{\CM}{\mathcal{M}}

 \def\CF{\mathcal F}
 \newcommand{\ZZ}{\mathbb{Z}}
\newcommand{\CR}{\mathcal R}
 \newcommand{\NN}{\mathbb N}

\newcommand{\co}{\circ}
\newcommand{\la}{\langle}
\newcommand{\ra}{\rangle}
 \gdef\tp{{\mathrm{tp}}}
 \gdef\CL{\mathcal{L}}
 \gdef\st{\operatorname{st}}
 \gdef\cl{\operatorname{cl}}
 \gdef\CO{\mathcal O}

 \gdef\CX{\mathcal{X}}
 \gdef\CY{\mathcal{Y}}

 \gdef\Rom{\RR_\mathrm{om}}
 \gdef\Lom{\CL_\mathrm{om}}
 \gdef\Rfull{\RR_\mathrm{full}}
 \gdef\Lfull{\CL_\mathrm{full}}

 \gdef\mfR{\mathfrak{R}}
 \gdef\fRfull{\mathfrak R_\mathrm{full}}
 \gdef\fRom{\mathfrak R_\mathrm{om}}
 \gdef\mfg{\mathfrak{g}}
 \gdef\gab{\mathfrak{g}_\mathrm{ab}}

\gdef\dcl{\operatorname{dcl}}
\gdef\ccdot{{\cdot}}
\gdef\SL{\mathscr{L}}
\gdef\SLm{\SL_{\mathrm{max}}}
\gdef\Gab{G_{\mathrm{ab}}}
\gdef\Gamab{\Gamma_{\mathrm{ab}}}

\gdef\piab{\pi_{\mathrm{ab}}}

 \newcommand{\bm}{{\mathbf m}}
\newcommand{\CS}{\mathcal S}

 \gdef\CP{\mathcal P}

 \gdef\bB{\overline{B}}
 \gdef\GL{\mathrm{GL}}
 \gdef\rhoab{\rho^{\mathrm{ab}}}
 \gdef\CFab{\CF^{\mathrm{ab}}}

\makeatletter \def\@tocline#1#2#3#4#5#6#7{\relax \ifnum
  #1>\c@tocdepth 
  \else
  \par \addpenalty\@secpenalty\addvspace{#2}%
  \begingroup \hyphenpenalty\@M \@ifempty{#4}{%
    \@tempdima\csname r@tocindent\number#1\endcsname\relax }{%
    \@tempdima#4\relax }%
  \parindent\z@ \leftskip#3\relax \advance\leftskip\@tempdima\relax
  \rightskip\@pnumwidth plus4em \parfillskip-\@pnumwidth
  #5\leavevmode\hskip-\@tempdima \ifcase #1 \or\or \hskip 1em \or
  \hskip 2em \else \hskip 3em \fi%
  #6\nobreak\relax \dotfill\hbox to\@pnumwidth{\@tocpagenum{#7}}\par
  \nobreak \endgroup \fi} \makeatother

\usepackage[sorted]{amsrefs}

\newtheorem*{thmnon}{Theorem}
\newtheorem{thm}{Theorem}[section]
\newtheorem{cor}[thm]{Corollary}
\newtheorem{lem}[thm]{Lemma}
\newtheorem{prop}[thm]{Proposition}
\newtheorem{claim}[thm]{Claim}
\newtheorem{fact}[thm]{Fact}

\theoremstyle{definition}
\newtheorem{defn}[thm]{Definition}
\newtheorem{ntn}[thm]{Notation}
\newtheorem{rem}[thm]{Remark}
\newtheorem{sample}[thm]{Example}

\numberwithin{equation}{section}

\usepackage[colorlinks=true, linkcolor=blue]{hyperref}

\title{Limits of definable families and dilations in nilmanifolds}
\author[Y.Peterzil]{Ya'acov Peterzil} \address{University of Haifa}
\email{kobi@math.haifa.ac.il} \author[S.Starchenko]{Sergei Starchenko}
\address{University of Notre Dame} \email{sstarche@nd.edu}
\reversemarginpar
\thanks{The first author was partially supported by ISF grant 290/19}
\thanks{The second author was partially supported by NSF research grant DMS-1800806}

\begin{document}
\maketitle

\begin{abstract} Let $G$ be a unipotent group and
  $\CF=\{F_t:t\in (0,\infty)\}$ a family of subsets of $G$, with $\CF$
  definable in an o-minimal expansion of the real field.  Given a
  lattice $\Gamma\subseteq G$, we study the possible Hausdorff limits of
  $\pi(\CF)$ in $G/\Gamma$ as $t$ tends to $\infty$ (here
  $\pi:G\to G/\Gamma$ is the canonical projection). Towards a
  solution, we associate to $\CF$ finitely many real algebraic
  subgroups $L\subseteq G$, and, uniformly in $\Gamma$, determine if the
  only Hausdorff limit at $\infty$ is $G/\Gamma$, depending on whether
  $L^\Gamma=G$ or not. The special case of polynomial dilations of a
  definable set is treated in detalis.

\end{abstract}
\section{Introduction}

\tableofcontents

Let $G$ be $\la \mathbb R^n,+\ra$ or more generally a real unipotent
group, and let $X\subseteq G$ be a definable set in some o-minimal
structure over $\mathbb R$. In \cite{o-minflows} and \cite{nilpotent}
we examined the following problem: \emph{For a lattice
  $\Gamma\subseteq G$, and $\pi:G\to G/\Gamma$, what is the topological
  closure of $\pi(X)$ in $G/\Gamma$?}

Using model theoretic machinery, we described the frontier of
$\cl(\pi(X))$ as the projection of finitely many definable families of
cosets of positive dimensional subgroups associated to $X$. The answer
can be seen, in a certain sense, as uniform in $\Gamma$.

Here we consider an extension of the problem:

\emph{For $G$ as above, let $\{F_s:s\in S\}$ be a family of subsets of
  $G$ that is definable in an o-minimal structure over the reals, let
  $\Gamma\subseteq G$ be a lattice and $\pi\colon G\to G/\Gamma$ the
  projection. What are the possible Hausdorff limits of the family
  $\{\pi(F_s):s\in S\}$ in $G/\Gamma$? How does the answer vary with
  $\Gamma$?}

Some results of this paper can be seen as an extension of work
\cite{KSS} and \cite{fish} on polynomial dilations in nilmanifolds.
But instead of considering equidistributionf of certain measures
linked to these dilations, we focus here on topological properties
(see Section \ref{sec: polynomial} below).

Our precise setting is as follows: Let $\Rom$ be an o-minimal
expansion of the field of reals.  Let $\CF=\{F_t:t\in (0,\infty)\}$ be
an $\Rom$-definable family of subsets of $G$, and let $\Gamma$ be a
lattice in $G$. We study the possible Hausdorff limits of the family
$\{\pi(F_t)\colon t\in (0,\infty)\}$, as $t$ tends to $\infty$.  Using
model theory, we replace the Hausdroff limits question by a question
on non-standard members of the family in an elementary extension. More
precisely, we consider an elementary extension $\mfR$ of
$\la \Rom,\Gamma\ra $ where for every definable set $Z$ in
$\la \RR,\Gamma\ra$ we denote by $Z^\sharp$ its realization in $\mfR$
(see Section \ref{sec: preliminaries} for details). Now every
Hausdorff limit at $\infty$ of $\pi(\CF)$ is the standard part of
$\pi(F_\tau^\sharp\ccdot \Gamma^\sharp)$ for $\tau >\RR$ a
non-standard parameter in $\mfR$ (see Section
\ref{sec:limits-at-infinity}). Thus, the problem reduces to the study
of sets of the form $\st(F_\tau^\sharp\ccdot \Gamma^\sharp)$.

Similarly to the answers to the closure problem, we associate to the
family $\CF$ finitely many normal co-commutative subgroups
$L_i\subseteq G$, and then for every $\Gamma$, the answers depend on
whether one of the $L_i$ is a $\Gamma$-dense subgroup or not. More
precisely, (see Section \ref{sec:prel-unip-groups} for details), let
$L^\Gamma$ be the smallest $\Gamma$-rational real algebraic subgroup
of $G$ containing $L$. We prove: (Theorem \ref{thm:main-hlin1}):

\begin{thmnon}[see Theorem~\ref{thm:main-hlin1}]

  Let $G$ be a unipotent group, $\CF=\{F_t:t\in (0,\infty)\}$ an
  $\Rom$-definable family of subsets of $G$.

  Then, there exists a finite collection $\SL(\CF)$ of normal
  co-commutative subgroups of $G$, such that for every lattice
  $\Gamma\subseteq G$ and $\pi:G\to G/\Gamma$, we have:

  \begin{enumerate}
  \item $L^\Gamma=G$ for some $L\in \SL(\CF)$ if and only if
    $\pi(\CF)$ converges strongly to $G/\Gamma$ at $\infty$
    (i.e. $G/\Gamma$ is the only Hausdroff limit at $\infty$ of
    $\pi(\CF)$ and this remains true for every lattice in $G$
    commensurable with $\Gamma$).

  \item $L^\Gamma\neq G$ for all $L\in \SL(\CF)$ if and only if there
    exists a subgroup $\Gamma_0\subseteq \Gamma$ of finite index such
    that all Hausdorff limits at $\infty$ of $\pi_0(\CF)$ are proper
    subsets of $G/\Gamma_0$ (here $\pi_0:G\to G/\Gamma_0$ is the
    quotient map).
  \end{enumerate}
\end{thmnon}

Note that the theorem above does not identify all the possible
Hausdorff limits of families $\pi(\CF)$ in $G/\Gamma$. However, we can
do it when $G$ is a abelian and $\CF$ is a family of polynomial
dilations with no constant term (see \ref{sec: polynomial}). We prove:

\begin{thmnon} [see Corollary \ref{cor:dil-def}]
  Let $\{ \rho_t\colon \RR^k\to \RR^m\colon t\in (0,\infty)\}$ be a
  family of polynomial dilations with no constant term, and
  $X\subseteq \RR^k$ an $\Rom$-definable set.

  Then there are linear subspaces $L_1,\ldots,L_n\subseteq \RR^n$, and
  there is a coset of a linear space $\bar c+V\subseteq (\RR^m)^n$ such
  that set of Hausdorff limits at $\infty$ of the family
  $\{\pi_\Gamma{\co}\rho_t(X) \colon t\in (0,\infty)\}$ is exactly the
  family
 \[\left\{\pi_\Gamma\Bigl ( \bigcup_{i=1}^n  (d_i+L_i^\Gamma)\Bigr ) : (d_1,\ldots, d_n)\in  \bar c+V^{\Gamma^n}\right\}.\]

 In particular, it is the projection under $\pi_\Gamma$ of a definable
 family of subsets of $\RR^n$.

\end{thmnon}

Our work on dilations was motivated by \cite{KSS}, of Kra, Shah and
Sun.

\medskip

\noindent{\bf The structure of the paper.} From a model theoretic
point of view, the main complexity of this work over the closure
theorems in \cite{nilpotent} is the fact that we study sets defined
over $\RR\la \tau\ra$, where $\tau$ is a non-standard parameter as
above.  This requires several adjustments to our previous work in
\cite{o-minflows} and \cite{nilpotent}.  In Section 2 we develop the
notion of short and long types (which replace bounded and unbounded
types over $\RR$).  In addition, we modify the theory of
$\mu$-stabilizers developed in \cite{mustab}, so it fits our setting.
In Section 3 we study types and their nearest co-commutative subgroups
(again, the results in \cite{nilpotent} need adjustments since the
types are over $\RR\la \tau\ra$). In Section 4, lattices come in and
we prove the main theorems about the $\Gamma$-closure of types.  In
Section 5, we study definable sets over $\RR\la \tau\ra$ and formulate
conditions under which such sets are $\Gamma$-dense.  In Section 6 we
translate the results obtained thus far back to the original problem
of Hausdorff limits, and in Section 7 we study in more details
families given by polynomial dilations.

\medskip

\noindent{\bf Acknowledgements.} We thank Amos Nevo for suggesting
this problem and explaining its ergodic theoretic origins.

\section{Long types and $\mu$-stabilizers}
\subsection{Model theoretic preliminaries}\label{sec: preliminaries}
As background on model theory and o-minimality we refer to \cite{omin}
and \cite{marker}. We follow the set-up from \cite
{o-minflows}*{Section 2} and \cite{mustab}*{Section 2.3}.

We fix an o-minimal structure $\Rom=\la \RR,<,+,\cdot,\cdots\ra$
expanding the real field and denote by $\Lom$ its language. For
convenience we add to $\Lom$ a constant symbol for every real number.

We use $\Lfull$ for a language in which every subset of $\RR^n$,
$n\in \mathbb N$, has a predicate symbol, and denote the corresponding
structure on $\RR$ by $\Rfull$. This will allow us to talk about
lattices as definable sets.

We let $\fRfull=\la \mfR,<\ldots\ra$ be an elementary extension of
$\Rfull$ which is $|\RR|^+$-saturated and
strongly-$|\RR|^+$-homogeneous, and let $\fRom$ be the reduct to
$\Lom$. Clearly, $\fRom$ is an elementary extensions of $\Rom$.

We use Roman letters $X,Y,Z$ to denote subsets of $\RR^n$ and let
$X^\sharp,Y^\sharp,Z^\sharp$ denote their realizations in
$\fRfull$. We use script $\CX$ to denote subsets of $\mfR^n$ which are
not necessarily of the form $X^\sharp$. When we write $A\subseteq \mfR$,
for a parameter set over which definable sets and types are
considered, we mean that $|A|\leq |\RR|$.

For $\CL=\Lom$ or $\CL=\Lfull$, as usual, a complete $\CL$-type over
$A$ is an ultrafilter on sets which are $\CL$-definable using
parameters in $A$.  For $A\subseteq \mfR$ and $\CX\subseteq \mfR^n$ an
$\Rom$-definable set over $A$, we let $S_\CX(A)$ be the collection of
all complete $\Lom$-types over $A$, containing the set $\CX$.  If
$\CX=X^\sharp$ for some $\Lom$ definable $X\subseteq \RR^n$ then
instead of $S_{X^\sharp}(A)$ we write $S_X(A)$.  For $p\in S_\CX(A)$
we let $p(\mfR)$ denote the set of its realizations in $\fRom$.

Unless otherwise stated, by ``definable'' we mean
``$\Lom$-definable''. In particular $\dcl$ denotes the definable
closure in the structure $\fRom$. Note that by our assumptions, $\Lom$
contains constant symbols for real numbers, hence, by definability of
Skolem functions, for any set $A\subseteq \mfR$, the definable closure
$\dcl(A)$ is an elementary substructure of $\fRom$ which contains
$\Rom$ as an elementary substructure.

A \emph{type-definable subset of $\mfR^n$, over $A$,} is the
intersection of (possibly infinitely many) definable sets over $A$,
which by our convention means $\Lom$-definable sets. The notion of a
$\Lfull$ type-definable set is similarly defined.  Since
$|A| \leq |\mathbb R|$, every collection of such definable sets is
bounded in size. A subset of $\mfR^n$ is said to be \emph{
  type-definable} if it is type-definable over some $A\subseteq \mfR$.

We let $\CO$ be the convex hull of $\RR$ in $\mfR$,
namely,
\[\CO=\{\alpha\in \mfR: \exists r\in \RR^{>0}\, |\alpha|<r\},\]
It is a valuation ring of $\mfR$, whose associated maximal ideal is
\[\bm=\{\alpha\in \mfR: \forall r\in \RR^{>0}\, |\alpha|<r\}.\]

The ring homomorphism $\CO\to \CO/\bm$ restricts to an isomorphism
between $\RR$ and $\CO/\bm$.  The corresponding ring homomorphism
$\st\colon \CO\to \RR \simeq\CO/\bm$ is called \emph{the standard part
  map}, and we extend it coordinate-wise to $\st:\CO^n\to \RR^n$. For
$\CX\subseteq \mfR^n$, we write $\st(\CX)$ instead of $\st(\CX\cap \CO^n)$.

For $a\in \RR^n$ and $r>0$, we let $B_r(a)=\{x\in \RR^n:|x-a|<r\}$.

We need the following lemma.

\begin{lem}\label{lem:st-closed}

  \begin{enumerate}
  \item If $\CX\subseteq {\mfR}^n$ is an $\Lfull$ type-definable set
    then $\st(\CX)$ is a closed subset of $\RR^n$.
  \item For a set $X\subseteq \RR^n$, we have $\cl(X)=\st(X^\sharp)$
    (where $\cl(X)$ is the topological closure of $X$).
  \item Let $\Sigma$ be a collection of $\Lfull$-definable subsets of
    ${\mfR}^n$ with $|\Sigma| \leq |\RR|$. If $\Sigma$ is closed under
    finite intersections, then
    \[ \st\Bigl(\bigcap \Sigma \Bigr) = \bigcap_{\CX\in \Sigma}
      \st(\CX). \] In particular $\st(\bigcap \Sigma)$ is closed.
  \end{enumerate}
\end{lem}
\begin{proof}

  (1) If $a\in\cl(\st(\CX))$ then for every $r\in \RR^{>0}$,
  $B_r(a)^\sharp\cap \st(\CX)\neq \emptyset$ and therefore also for
  every $r\in \RR^{>0}$, $B_r(a)^\sharp\cap \CX\neq \emptyset$. By
  saturation, there is $b\in \CX$ such that $b\in a+\bm$, so
  $a\in \st(\CX)$.

  (2) is easy and (3) is just \cite[Claim 3.1]{o-minflows}
\end{proof}

The following standard fact is easy to prove.
\begin{fact}\label{fact:st-com}
  Let $X\subseteq \RR^m$ and $Y\subseteq \RR^n$ be closed subsets and
  $f\colon X\to Y$ a continuous function. Let
  $\st_m\colon \CO^m\to \RR^m$ and $\st_n\colon \CO^n\to \RR^n$ be the
  corresponding standard part maps. For every
  $\gamma\in \CO^m\cap X^\sharp$ we have
  $f(\st_m(\gamma))=\st_n(f(\gamma))$.  In particular, for every
  $\CX \subseteq X^\sharp$ we have
  $f(\st_m(\CX))\subseteq \st_n(f(\CX))$.
\end{fact}

\medskip

We will need the following important result:
\begin{fact}[{\cite{lou-limit}*{Proposition 8.1}}]\label{fact:st-def}
  If $\CX$ is definable in $\fRom$ then $\st(\CX)$ is definable in
  $\Rom$.
\end{fact}

\medskip

Let $G\subseteq \RR^m$ be an $\Rom$-definable group, namely the universe of
$G$ and the group operation are definable in $\Rom$. For example, any
real algebraic, or more generally semi-algebraic group is definable in
$\Rom$ (for more on definable groups in o-minimal structures see
\cite{Otero}). The group $G$ can be endowed with a group topology with
a definable basis (see \cite{p}). This topology might disagree with
the natural o-minimal topology, coming from the fact that $G$ is a
subset of $\RR^m$. However, as we observed in \cite{mustab}*{Claim
  3.1}, we may embed $G$ definably as a closed subset of $\RR^n$ for
some $n$, such that the above group topology agrees with the induce
$\RR^n$-topology. Thus, whenever $G\subseteq \RR^n$ is an $\Rom$-definable
group we assume it to be closed in $\RR^n$, and in addition assume
that the Euclidean topology makes $G$ a topological group.

We will use very often the following.
\begin{fact}[\cite{p}]\label{fact:group-hom-cont}
  Let $G$ be an $\Rom$-definable group.
  \begin{enumerate}
  \item If $H\subseteq G$ is an $\Rom$-definable subgroup then $H$ is
    closed in $G$.
  \item If $H$ is an $\Rom$-definable group and $f\colon G\to H$ an
    $\Rom$-definable homomorphism then $f$ is continuous.
  \end{enumerate}
\end{fact}

For $G\subseteq \RR^n$ as above, we consider two distinguished subgroups of
$G^\sharp$. The first is the infinitesimal group $\mu_G$, defined as
follows:
\[\mu_G=\bigcap \{X^\sharp: X\subseteq G\mbox{ an $\Rom$-definable open
    neighborhood of $e$}\}.\] It is a type-definable subgroup and
under our assumption on $G$ it equals $e+\bm^n\cap G^\sharp$, with
$\bm\subseteq \CO\subseteq \mfR$ the infinitesimal ideal defined above.

The second subgroup is $\CO_G$, defined by:
\[ \CO_G=\bigcup\{X^\sharp: X\subseteq G \mbox{ an $\Rom$-definable compact
    neighborhood of $e$}\}.\]
$\CO_G$ is a $\bigvee$-definable (or
Ind-definable) subgroup of $G^\sharp$, which equals, under our
assumptions on $G$, to $\CO^n\cap G^\sharp$. In particular,
$G\subseteq \CO_G$. The group $\mu_G$ is a normal subgroup of $\CO_G$ and
the latter can be written as the semi-direct product
$\CO_G=\mu_G\rtimes G$. We identify the quotient $\CO_G/\mu_G$ with
$G$ and call the quotient map $\st_G:\CO_G\to G$ \emph{the standard
  part map}. As before, we extend it coordinate-wise to
$\st_G\colon \CO_G^n \to G^n$. By our assumptions on $G$, for every
$a\in \CO_G$, we have $\st(a)=\st_G(a)\in G$.  In particular, Lemma
\ref{lem:st-closed} holds if one restricts to subsets of $G$ and to
$\st_G$.

   When the underlying group $G$ is fixed we omit the subscript $G$
   and just use $\CO$, $\mu$ and $\st$.

   As before, given an arbitrary set $\CX\subseteq G^\sharp$, we let
   $\st(\CX)$ denote the set $\st(\CX\cap \CO_G)\subseteq G$.

   \begin{fact}\label{fact:maps} Assume that $G_1,G_2$ are definable in
     $\Rom$. Then,

     \begin{enumerate} \item
       $\mu_{G_1\times G_2}=\mu_{G_1}\times \mu_{G_2}$ and
       $\CO_{G_1\times G_2}=\CO_{G_1}\times \CO_{G_2}$.

     \item If $f\colon G_1\to G_2$ is an $\Rom$-definable surjective
       homomorphism then $f(\mu_{G_1})=\mu_{G_2}$ and
       $f(\CO_{G_1})=\CO_{G_2}$
     \end{enumerate}
   \end{fact} \proof (1) follows from the fact that the topology on
   $G_1\times G_2$ is the product topology.  For (2), see
   \cite{nilpotent}*{Lemma 3.10}.\qed

   \subsection{Long and short sets}

   Let $G\subseteq \RR^n$ be an $\Rom$-definable group.

   \begin{defn} A subset $\CX\subseteq G^\sharp$ is called \emph{left-short
       in $G$} if there exists some compact set $K\subseteq G$ and
     $g\in G^\sharp$ such that $\CX\subseteq K^\sharp\ccdot g$ (since $G$
     can be written as an increasing union of relatively compact
     $\Rom$-definable open sets, we may always take $K$ to be
     $\Rom$-definable).

     Otherwise, $\CX$ is called \emph{left-long in $G$}. For
     $A\subseteq \mfR$, we say that a type $p\in S_G(A)$ is left-short
     (left-long) in $G$ if $p(\mfR)$ is left-short (left-long) in $G$.
   \end{defn}

   The following are easy to verify:

   \begin{lem} \label{lem:properties of long} Given
     $\CX\subseteq G^\sharp$,
     \begin{enumerate}
     \item $\CX$ is left-short in $G$ if and only if
       $\CX\ccdot\CX^{-1}\subseteq K^\sharp$ for some compact set
       $K\subseteq G$.
     \item For every $g\in G^\sharp$, $\CX$ is left-short in $G$ if
       and only if $\CX g$ is left-short in $G$.
     \item For every $g\in G$, $\CX$ is left-short in $G$ if and only
       if $g\CX$ is left-short in $G$.
     \item If $\CX=X^\sharp$, for $X\subseteq G$ an $\Rfull$-definable set,
       then $\CX$ is left-short in $G$ if and only if $X$ is bounded
       in $\RR^n$.
     \end{enumerate}

   \end{lem}

   We may similarly define right-short and right-long in $G$ and in
   general these notions are different. {\bf However, for the rest of
     the paper we use \emph{short} and \emph{long} to refer only to
     left-short and left-long}.

   If $H\subseteq G$ is an $\Rom$-definable subgroup and $\CX$ a subset of
   $H^\sharp$ then, by Lemma~\ref{lem:properties of long}(1), $\CX$ is
   short in $H$ if and only if it is short in $G$, hence we omit the
   reference to the group when the context is clear.

   Note that by saturation, an $\Lfull$ type-definable set
   $\CX\subseteq G^\sharp$ is short if and only if $\CX\subseteq \CO \ccdot g $
   for some $g\in G^\sharp$.

   \begin{lem}\label{lem:short1} Let $G_1,G_2, G$ be $\Rom$-definable
     groups.
     \begin{enumerate}
     \item If $\CX\subseteq G_1^\sharp$ is short and $f:G_1\to G_2$ is an
       $\Rom$-definable homomorphism then $f(\CX)$ is short in $G_2$.

     \item If $\CX_1\subseteq G_1^\sharp$ is short and
       $\CX_2\subseteq G_2^\sharp$ is short then $\CX_1\times \CX_2$ is
       short in $G_1\times G_2$.

     \item If $H_1$, $H_2$ are two normal $\Rom$-definable subgroups
       of $G$ and $\CX\subseteq G^\sharp$ an arbitrary set then the image
       of $\CX$ in $G^\sharp/(H_1^\sharp\cap H_2^\sharp)$ is short if
       and only if its images in $G^\sharp/H_1^\sharp$ and in
       $G^\sharp/H_2^\sharp$ are short.
     \end{enumerate}
   \end{lem}
   \begin{proof} (1) and (2) are immediate since the image of a
     compact set under a $\Rom$-definable homomorphism is compact, and
     similarly the direct product of such sets is compact.

     For (3), let $\pi_i:G\to G/H_i$, $i=1,2$, be the natural
     projections, and let $\pi:G\to G/H_1\times G/H_2$ be the map
     $\pi(g)=(\pi_1(g),\pi_2(g))$.  The kernel of $\pi$ is
     $H_1\cap H_2$ hence the image is isomorphic to $G/(H_1\cap H_2)$.

     If $\pi_1(\CX)$ and $\pi_2(\CX)$ are both short then by (2), so
     is
     $\pi_1(\CX)\times \pi_2(\CX)\subseteq G^\sharp/H_1^\sharp\times
     G^\sharp/H_2^\sharp$. But $\pi(\CX)$ is contained in
     $\pi_1(\CX)\times \pi_2(\CX)$ so also short. The converse follows
     from (1) using the natural homomorphisms from $G/(H_1\cap H_2)$
     onto $G/H_i$, $i=1,2$.
   \end{proof}

   \begin{lem}\label{lem:short2} For $A\subseteq \mfR$, and
     $p\in S_G(A)$, we have:\begin{enumerate}

     \item $p$ is short in $G$ if and only if there exists
       $a\in \dcl(A)\cap G^\sharp$ such that $p\vdash \mu \ccdot a$,
       namely $p(\mfR)\subseteq \mu\ccdot a$.

     \item Let $H\subseteq G$ be a $\Rom$-definable normal subgroup and
       $\pi:G\to G/H$ the quotient map. Then $\pi(p)\in S_{G/H}(A)$ is
       short in $G/H$ if and only if there exists
       $a\in \dcl(A)\cap G^\sharp$ such that
       $p(\mfR)\subseteq \mu \ccdot a H^\sharp$.
     \end{enumerate}
   \end{lem}
   \proof (1) Assume that $p$ is short. Hence
   $p(\mfR)\ccdot p(\mfR)^{-1}\subseteq K^\sharp$ for some $\Rom$-definable
   compact set $K\subseteq G$.  By logical compactness, there exists an
   $A$-definable set $\CX$ in $p$ such that
   $\CX \ccdot\CX^{-1}\subseteq K$. By definability of Skolem functions in
   o-minimal structures, the set $\CX$ contains a point
   $a\in \dcl(A)$. Consider the complete $\Rom$-type over $A$,
   $p{\ccdot} a^{-1}$. We have
   $p(\mfR)\ccdot a^{-1} \subseteq \CX\ccdot a^{-1}\subseteq
   K^\sharp$. Let $\beta\models p\ccdot a^{-1}$ and $g=\st(\beta)$
   (this is defined since $\beta\in \CO$). Because $p\ccdot a^{-1}$ is
   a complete $A$-type and $\beta\in \mu\ccdot g$, we have
   $p\ccdot a^{-1}\vdash \mu{\ccdot}g$. Hence
   $p(\mfR) \subseteq \mu \ccdot g \ccdot a$. Clearly,
   $g\ccdot a \in \dcl(A)$.

   The converse is clear.

   For (2), notice that if $p$ is short then by
   Lemma~\ref{lem:short1}(1), $\pi(p)$ is short and hence by (1) there
   exists $b\in \dcl(A)\cap (G/H)^\sharp$ such that
   $\pi(p(\mfR))\in \mu_{G/H}\ccdot b$.

   We now take any $a\in\dcl(A)$ in the $A$-definable set
   $ \pi^{-1}(b)$, and we have
   $p(\mfR)\subseteq \pi^{-1}(\mu_{G/H}\ccdot b) \subseteq \mu_G\ccdot a
   H^\sharp$.

   For the converse, notice that by Fact~\ref{fact:maps},
   $\pi(\mu_G)=\mu_{G/H}$.  Hence,
   $\pi(\mu_G\ccdot a)=\mu_{G/H}\ccdot \pi(a)$, so if
   $p(\mfR)\subseteq \mu_G\ccdot aH^\sharp$ then
   $\pi(p)(\mfR) \subseteq \mu_{G/H}\ccdot \pi(a)$ is short.  \qed

   \subsection{The $\mu$-stabilizer of a type}

   We fix an $\Rom$-definable group $G$.  In \cite{mustab} we
   developed a theory for $\mu$-stabilizers of types over $\RR$.  Here
   we take a more general viewpoint which we now explain.

   Consider the set $S_G(A)$. Given $p\in S_G(A)$ we let
   $\mu \ccdot p$ (below written as $\mu p$) denote the partial type
   over $A$ whose realization is the set $\mu \ccdot p(\mfR)$.

   \begin{rem} We note that when we consider here types over arbitrary
     $A\subseteq \mfR$, then, unlike \cite{mustab}, we still keep
     $\mu=\mu_G$ fixed and not change it to a smaller infinitesimal
     group (namely, the intersection of all $A$-definable open
     neighborhoods of $e$).

     Since our point of view here is slightly different from
     \cite{mustab}, we go briefly through the results we need and
     explain how their proofs differ from the analogous results in
     \cite{mustab}.
   \end{rem}

   Given $p,q\in S_G(A)$, we say that $p$ and $q$ are
   $\mu$-equivalent, $p\sim_\mu q$, if $\mu p=\mu q$, i.e.
   $\mu p (\mfR)=\mu q (\mfR)$.
   \begin{fact}\label{fact:mutypes} For $p,q\in S_G(A)$, the
     following are equivalent:

     \begin{enumerate}
     \item $p\sim_\mu q$
     \item $\mu p(\mfR) \cap \mu q(\mfR)\neq \emptyset$.
     \end{enumerate}
   \end{fact}
   \noindent(see \cite{mustab}*{Claim 2.7} for an identical argument).

   It is easy to verify that if $p\sim_\mu q$ then $p$ is a long if
   and only if $q$ is long.

   Let $S^\mu_G(A)=\{\mu p:p\in S_G(A)\}$. The group $G$ acts from the
   left on $S^\mu_G(A)$ by $g\ccdot \mu p=\mu (g p)$.  The following
   subgroup plays a crucial role in our analysis:
   \begin{defn}
     Given $p\in S_G(A)$, the \emph{left stabilizer of $\mu p$} is
     defined as:
\[ \Stab^\mu(p)=\{g\in G: g\ccdot \mu p=\mu p\}.\]

\end{defn}

Since the definition of $\Smu(p)$ depends only on $\mu p$, if
$p\sim_\mu q$ then $\Smu(p)=\Smu(q)$.

Our main focus in \cite{mustab} was on unbounded definable types.  For
$\Rom$-definable groups, and $\Lom$-types over $\RR$ these are types
which do not contain any formula over $\RR$ defining a compact subset
of $G$. Since we are considering here types which are not only over
$\RR$ our focus is shifted to long types.

Recall that for $p$ an $\Lom$-type we let $\dim(p)$ be the smallest
o-minimal dimension of the formulas in $p$.
\begin{defn} We say that a type $p\in S_G(A)$ is \emph{$\mu$-reduced}
  if for all $q\in S_G(A)$, if $p\sim_\mu q$ then
  $\dim(p)\leq \dim(q)$.\end{defn}

Clearly, every $p\in S_G(A)$ is $\mu$-equivalent to a $\mu$-reduced
type in $S_G(A)$: just take a $\mu$-equivalent type of minimal
dimension (but there might be more than one such). Notice that by
Lemma~\ref{lem:short2}, if $p$ is short and $\mu$-reduced then
$\dim p=0$ and $p=\tp(a/A)$ for some $a\in\dcl(A)$.

Our main goal in this section is to prove:

\begin{prop}\label{prop:main-mustab} Let $p\in S_G(A)$.
  Then \begin{enumerate} \item $\Smu(p)$ is $\Rom$-definable and can
    be written as $\st(\CS\ccdot\alpha^{-1})$, for some definable
    $\CS$ in $p$ and $\alpha\models p$.  \item If $p$ is a long type
    then $\dim(\Smu(p))>0$. Moreover, in this case $\Smu(p)$ is a
    torsion-free solvable group.

  \end{enumerate}

\end{prop}

The proof is very similar to the proof of \cite{mustab}*{Theorem 3.10}
so we only point out the differences. As we noted above, we may assume
that $p$ is $\mu$-reduced and if $p$ is short that
$\mu p =\mu \ccdot a$ for some $a\in\dcl(A)$ so $\Smu(p)$ is trivial.
Thus, we fix a long $\mu$-reduced type $p\in S_G(A)$ and
$\alpha\in p(\mfR)$.

We start with an analogue of \cite{mustab}*{Claim 3.12}:

\begin{claim}\label{claim:like 3.12} If $Y\subseteq G^\sharp$ is $A$-definable and $\dim Y<\dim p$ then $\CO \ccdot\al\cap Y=\emptyset$.\end{claim}
\proof Assume towards contradiction that
$\beta\in Y\cap \CO\ccdot\alpha$. Then
$\beta \in \mu \ccdot r\alpha $, for some $r\in G$, hence
$r^{-1}\beta \in \mu p$.  But
$\dim(r^{-1}\beta/A)\leq \dim(Y)<\dim p$, contradicting the fact that
$p$ is $\mu$-reduced.\qed

Next, we note, just like \cite{mustab}*{Claim 3.8}, that for every
$A$-definable set $\CS$ in $p$, $\Smu(p)\subseteq \st(\CS\ccdot \al^{-1})$.
Indeed, if $g\in \Smu(p)$ then there exists $\beta\models p$ and
$\eps\in \mu$ such that $g\al=\eps\beta$. It follows that
$\beta\in \CS$ and $\beta\al^{-1}\in \CO$, thus
$g=\st(\beta \al^{-1})\in \st(\CS\ccdot \al^{-1})$.

The next claim is similar to \cite{mustab}*{Claim 3.13}.

\begin{claim}\label{claim:like 3.13}  There exists an $A$-definable set $\CS$ in $p$ such that every element in $\CS\cap \CO\ccdot \al$ realizes $p$.\end{claim}
Let us explain the proof: As in \cite{mustab}, for every $A$-definable
set $\CS$ in $p$, the set $\CS\ccdot \al^{-1}\cap \CO$ is a relatively
definable subset of $\CO$. Hence, by \cite{mustab}*{Theorem B.2} it
has finitely many connected components (see precise definition of
connectedness there). We choose an $A$-definable such cell $\CS$ in
$p$ with $\dim S=\dim (p)$, for which the number of components of
$\CS\ccdot \al^{-1}\cap \CO$ is minimal. Using Claim~\ref{claim:like
  3.12}, we can prove, just as in \cite{mustab}*{Claim 3.13}, that any
$\beta\in \CS\cap \CO\ccdot \al$ must realize $p$.

Finally, we prove an analogue of \cite{mustab}*{Claim 3.14}:
\begin{claim}\label{claim:like 3.14} For $\CS$ as in
  Claim~\ref{claim:like 3.13}, we have
  $\Smu(p)=\st(\CS\ccdot \al^{-1})$.
\end{claim}
\proof It is sufficient to show that
$\st(\CS\ccdot\al^{-1})\subseteq \Smu(p)$, so we take
$g\in \st(\CS\ccdot\al^{-1})$ and note that for some $\eps\in \mu$, we
have $\eps g\al\in \CS\cap \CO\ccdot \al$, so by our choice of $\CS$,
$\eps g\al\models p$. It follows that
$\mu gp(\mfR)\cap p(\mfR)\neq \emptyset$, so by
Fact~\ref{fact:mutypes}, $g\ccdot \mu p=\mu p$. \qed

Thus, by Fact~\ref{fact:st-def}, $\Smu(p)=\st(\CS\ccdot \alpha^{-1})$
is definable.

Since $p$ is long the set $\CS\ccdot \al^{-1}$ is not contained in
$\CO$, thus $\Smu(p)$ is unbounded in $G$.

To see that $\Smu(p)$ is solvable, torsion-free we repeat the argument
from \cite{mustab}*{Theorem 3.6}: By \cite{mustab}*{Fact 3.25}, $G$
can be written as a product of two sets $G=C\ccdot H$, with $C\subseteq G$
a $\Lom$-definable compact set and $H$ a $\Lom$-definable torsion-free
solvable group. Thus, $\al\in G^\sharp$ as above can be written as
$\alpha=\eps\ccdot g\ccdot h^*$ for
$\eps\in \mu_G, g\in C, h^*\in H^\sharp$. It follows that
$\al\in \mu \ccdot(H^g)^\sharp \ccdot g$, so $\tp(\al/A)$ is
$\mu$-equivalent to a type $q\vdash (H^g)^\sharp \ccdot g$. But then
$\Smu(p)=\Smu(pg)\subseteq H^g$ so $\Smu(p)$ is a torsion-free solvable
group.

This ends the proof of Proposition~\ref{prop:main-mustab}.\qed

\begin{rem} In fact, the remainder of the proof of
  \cite{mustab}*{Theorem 3.12} goes through identically and thus we
  could have proved the stronger result, saying that for $p$ a
  $\mu$-reduced type over $A$, the dimension of $\Smu(p)$ equals to
  $\dim(p)$. However, this will not be needed here.
\end{rem}

\section{Nearest cosets}\label{sec:nearest-cosets}

We now assume again that $G$ is a definable group in $\Rom$.

\subsection{Nearest co-commutative cosets}\label{sec:near-co-comm}

\begin{defn}
  Given a type $p\in S_G(A)$, an $\Rom$-definable subgroup $H\subseteq G$
  and $a\in\dcl(A)\cap G^\sharp$, we say that the coset $aH^\sharp$
  \emph{is near $p$} if $p(\mfR)\subseteq \mu \ccdot a H^\sharp$.

  Sometimes we omit $\sharp$, write $p\vdash \mu a H$, and say that
  $aH$ is near $p$.
\end{defn}

Notice that in the above definition the subgroup $H$ is defined over
$\RR$, but the element $a$ is taken from $\dcl(A)\subseteq \mfR$.

\begin{rem}\label{rem:nearshort}
  By Lemma~\ref{lem:short2}, a type $p\in S_G(A)$ is short if and only
  if, for the trivial subgroup $\{e\}$, a coset $a\ccdot e$, is near
  $p$.

  Also, for a normal $\Rom$-definable subgroup $H\subseteq G$, some
  coset $aH$ is near $p$ if and only if the image of $p$ in $G/H$ is
  short.
\end{rem}

\begin{lem}\label{lem:near-inter}
  Let $p\in S_G(A)$, $H_1,H_2 \subseteq G$ be two $\Rom$-definable
  normal subgroups, $a_1,a_2\in \dcl(A)$, and assume both $aH_1$ and
  $aH_2$ are near $p$.  Then there exists $d\in\dcl(A)$ such that the
  coset $d(H_1\cap H_2)$ is near $p$.
\end{lem}
\proof Let $G_i=G/H_i$, $i=1,2$, and $\pi_i\colon G\to G_i$ the
natural projection. Let $f\colon G\to G_1\times G_2$ be the definable
homomorphism $f(g)=(\pi_1(g),\pi_2(g))$.  We have
$\ker(f)=H_1\cap H_2$.

By Lemma~\ref{lem:short2}(2), the images of $p(\mfR)$ in both $G/H_1$
and in $G/H_2$ are short. Hence, by Lemma~\ref{lem:short1}(3), its
image in $G/(H_1\cap H_2)$ is also short. By Lemma~\ref{lem:short2}
(2), there exists $d\in\dcl(A)$ such that $d(H_1\cap H_2)$ is near
$p$. \qed

In the case of a unipotent group $G$ and $A=\RR$, the above lemma
holds without assuming normality of $H_1$ and $H_2$ (see
\cite{nilpotent}*{Theorem 3.7}).  Unfortunately, in general, this
fails for an arbitrary $A$:

  \begin{sample}
    We consider the Heisenberg group, identified with $\RR^3$, as
    \[ [a,b,c]\ccdot [d,e,f]=[a+d,b+e,ae+c+f].\] We let
    $H_1=Z(G)=\{[0,0,x]: x\in \RR\}$, and $H_2=\{[0,t,t]: t\in
    \RR\}$. We now consider $H_1^\sharp$ and $H_2^\sharp$ in
    $G^\sharp$.  Fix $\tau \in \mfR$ with $\tau>\RR$ and let
    $A=\dcl(\tau)$.

    Consider the $1$-type over $\dcl(A)$:
\[q(t)=\{r<t<c: r\in \RR\, , c\in\dcl(A) \mbox{ with } c> \RR\}.\]
  Let $p(t)$ be the type over $\dcl(A)$ given by
  $\{[\tau,0,t]: t\models q\}$. For $\alpha=[\tau,0,0]$, the
  realizations of $p$ are contained in the coset $\alpha H_1$. It is
  easy to see that $p$ is a long type and we claim that $\alpha H_2$
  is near $p$.

  Indeed, consider the type $q_0=(1/\tau) q$. It is also a $1$-type
  over $\dcl(A)$, whose realizations are contained in $\mu\subseteq \mfR$,
  and, for every $\beta\models q_0$, the element
  $g_\beta=[0,\beta,\beta]$ is in $H_2^\sharp$.  Now, for
  $\varepsilon_\beta=[0,-\beta,-\beta]\in \mu_G$, we have
\[\varepsilon_\beta\ccdot \alpha\ccdot g_\beta= [0,-\beta,-\beta]\ccdot
  [\tau, \beta,\tau \beta+\beta]=[\tau,0,\tau \beta]\models p,\]
and
also
$\varepsilon_\beta\ccdot \alpha\ccdot g_\beta\in \mu_G \ccdot
\alpha\ccdot H_2^\sharp$.  Hence the coset $\alpha H_2$ is near $p$.

Thus both $\alpha H_1$ and $\alpha H_2$ are near $p$. However, since
$p$ is long, a coset of $H_1\cap H_2=\{e\}$ can not be near $p$.
\end{sample}

The main part of this paper deals with unipotent groups, and, in the
unipotent case, instead of nearest cosets, as in \cite{nilpotent}, it
is more convenient to work with nearest co-commutative cosets.

\begin{defn}
  We say that a subgroup $H\subseteq G$ is \emph{co-commutative} if it
  is normal and the quotient $G/H$ is abelian (equivalently $H$
  contains $[G,G]$).
\end{defn}

Since the intersection of two co-commutative subgroups is
co-commutative, using Lemma~\ref{lem:near-inter}, we may conclude:
\begin{cor}\label{cor:nearest-cor} Given $p\in S_G(A)$ there exists a
  smallest (by inclusion) $\Rom$-definable co-commutative subgroup
  $L\subseteq G$ such that for some $a\in \dcl(A)$ the coset $aL$ is near
  $p$.
\end{cor}

We can now define:
\begin{defn}
  Given $p\in S_G(A)$, \emph{a nearest co-commutative coset to $p$} is
  a coset of the form $a L$, where $a\in\dcl(A)$ and $L\subseteq G$ is an
  $\Rom$-definable co-commutative subgroup as in
  Corollary~\ref{cor:nearest-cor}.  It is unique up to $\mu_G$, namely
  if $a_1 L_1$ and $a_1 L_1$, are both nearest co-commutative cosets
  to $p$ then $L_1=L_2$ and
  $\mu_G\ccdot a_1 L_1^\sharp=\mu_G\ccdot a_2 L_2^\sharp$.

  We will denote this subgroup $L$ as $L_p$.
\end{defn}

We now prove some basic properties of nearest co-commutative cosets.

\begin{lem}\label{lem:properties} Let $p\in S_G(A)$.
  If $K\subseteq G$ is a compact definable set, $a\in\dcl(A)$ and
  $p(\mfR)\subseteq K^\sharp \ccdot a \ccdot L^\sharp$ for some
  $\Rom$-definable co-commutative $ L\subseteq G$ then $L_p\subseteq L$.
\end{lem}
\begin{proof}
  Clearly $\pi(p)$ is short in $G/L$, where $\pi\colon G \to G/L$ is
  the quotient map.  By Lemma~\ref{lem:short2}(2), there exists
  $a'\in\dcl (A)$ such that $a'\ccdot L$ is near $p$, hence
  $L_p\subseteq L$.
\end{proof}

\begin{lem}\label{lem:near-hom} Let $f\colon G\to H$ be an
  $\Rom$-definable surjective homomorphism of definable groups and
  $A\subseteq \mfR$. For a type $p\in S_G(A)$ and $q=f(p)$, if
  $D_p=aL_p$ is a nearest co-commutative coset to $p$ then $f(D_p)$ is
  a nearest co-commutative coset to $q$, and in particular,
  $L_q=f(L_p)$.
\end{lem}
\begin{proof} Since $f$ is surjective, it maps a co-commutative
  subgroup onto a co-commutative subgroup.

  Let $D_q$ be a nearest co-commutative coset to $q$.  It is
  sufficient to see that $\mu_HD_q=\mu_Hf(D_p)$.  We have
  $p\vdash \mu_G D_p$, so by Fact~\ref{fact:maps},
  $q\vdash \mu_H f(D_p)$, hence $D_q\subseteq \mu_H f(D_p)$. Conversely,
  since $q\vdash \mu_H D_q$ then $p\vdash \mu_G f^{-1}(D_q)$, so
  $D_p\subseteq \mu_G f^{-1}(D_q)$, hence $f(D_p)\subseteq \mu_H D_q$.
\end{proof}

\begin{lem}\label{lem:mu-nearest} For $p\in S_G(A)$, let $H\subseteq G$ be the
  $\mu$-stabilizer of $p$, and let $aL_p$ be a nearest co-commutative
  coset to $p$. Then $H\subseteq L_p$.
\end{lem}
\proof Fix $\beta\models p$. Then by assumption, there exists
$\varepsilon\in \mu$ and $\ell\in L_p^\sharp$ such that
$\beta=\varepsilon a \ell$.  Given $h\in H$, we have
$h \beta\in \mu p(\mfR)$, hence $h \beta=\varepsilon' a \ell'$, with
$\varepsilon'\in \mu$ and $\ell'\in L_p^\sharp$.  Thus
\[ h = h \beta \beta^{-1} = \varepsilon' a \ell' \ell^{-1} a^{-1}
  \varepsilon^{-1}.\] Since $L_p$ is normal in $G$ (so $L_p^\sharp$
normal in $G^\sharp$), it follows that $h\in \mu\ccdot
L_p^\sharp$. However, $h$ is in $G$ and $L_p$ is closed in $G$,
therefore $h\in L_p$.\qed

\subsection{The set $\SLm(\CX)$}\label{sec:set-clmax}

Again we fix a group $G$ definable in $\Rom$.

Recall that by our assumption, $G\subseteq \RR^n$ is a closed subset.
For $r>0$, we will denote by $\bB_r\subseteq G$ the set
$\bB_r(e)\cap G$, where $\bB_r(e)$ is the closed ball of radius $r$
centered at $e$. Clearly, each $\bB_r$ is a compact subset of $G$,
definable in $\Rom$, with
$\displaystyle \mu_G =\bigcap_{r\in \RR^{>0}} \bB_r^\sharp$.

\begin{lem}\label{lem:SX-ext}
  For $A\subseteq B \subseteq \mfR$, let $\CX\subseteq G$ be a set
  $\Lom$-definable over $A$, and $p\in S_\CX(A)$.
  \begin{enumerate}
  \item If $q\in S_\CX(B)$ is an extension of $p$ then
    $L_q \subseteq L_p$.
  \item There is $q\in S_\CX(B)$ extending $p$ such that $L_q=L_p$.
  \end{enumerate}
\end{lem}
\begin{proof}
  (1).  Choose $a_p\in \dcl(A)$ such that
  $p(\mfR) \subseteq \mu\ccdot a_p \ccdot L_p^\sharp$. We have
  $q(\mfR) \subseteq p(\mfR) \subseteq \mu\ccdot a_p \ccdot
  L_p^\sharp$. Hence the coset $a_pL_p$ is near $q$ and
  $L_q \subseteq L_p$.

  (2).  Let $\mathcal{Q}$ be the set of all $q\in S_\CX(B)$ extending
  $p$.  For each $q\in \mathcal{Q}$ we choose $b_q\in \dcl(B)$ such
  that $q(\mfR) \subseteq \mu \ccdot b_q \ccdot L_q^\sharp$.  We have
  \[ p(\mfR) \subseteq \bigcup_{q\in \mathcal{Q}} q(\mfR) \subseteq
    \bigcup_{q\in \mathcal{Q}} \mu \ccdot b_q \ccdot L_q^\sharp
    \subseteq \bigcup_{q\in \mathcal{Q}} \bB_1^\sharp\ccdot b_q\ccdot
    L_q^\sharp. \] Thus the type definable set $p(\mfR)$ is covered by
  a bounded family of definable sets of the form
  $\bB_1^\sharp\ccdot b_q\ccdot L_q^\sharp$.  Hence, by logical
  compactness, we can find a set $\CX_0\in p$, definable over $A$, and
  a finite subset $\mathcal{Q}_0\subseteq \mathcal{Q}$ such that
  \[ \CX_0 \subseteq \bigcup_{q\in \mathcal{Q}_0} \bB_1^\sharp\ccdot
    b_q\ccdot L_q^\sharp. \] Since $\dcl(A)$ is an elementary
  substructure of $\fRom$, we can find $a_q\in \dcl(A)$, for each
  $q\in \mathcal{Q}_0$, such that
  $\CX_0 \subseteq \bigcup_{q\in \mathcal{Q}_0} \bB_1^\sharp\ccdot
  a_q\ccdot L_q^\sharp$.  Since $p$ is a complete over $A$, there is
  $q\in \mathcal{Q}_0$ with
  $p(\mfR) \subseteq \bB_1^\sharp\ccdot a_q\ccdot L_q^\sharp$. By
  Lemma~\ref{lem:properties}, $L_p \subseteq L_q$, hence by (1), we
  have $L_p=L_q$.
\end{proof}

For $A\subseteq \mfR$ and a set $\CX\subseteq \mfR^n$ definable over $A$,
we denote by $\SL_A(\CX)$ the set
\[ \SL_A(\CX) =\{ L_p \colon p\in S_\CX(A) \}. \]

\begin{cor}\label{cor:incl-L} For $A\subseteq B \subseteq\mfR$, let
  $\CX\subseteq G^\sharp$ be definable over $A$. Then,
  \begin{enumerate}
  \item $\SL_A(\CX) \subseteq \SL_B(\CX)$.
  \item An $\Rom$-definable co-commutative subgroup $L$ of $G$ is
    maximal (by inclusion) in $\SL_A(\CX)$ if an only if it is maximal
    in $\SL_B(\CX)$
  \end{enumerate}
  \begin{proof}
    Follows from Lemma~\ref{lem:SX-ext}.
  \end{proof}
\end{cor}
\begin{rem}
  In general, for $A\subseteq B$ we do not have equality of sets,
  $\SL_A(\CX) =\SL_B(\CX)$. As an example, consider the group
  $G=(\RR^2,+)$ with $\CX=\{ (x,y)\in \mfR^2 \colon x\geq 0, y=x^2\}$.
  For $A=\RR$ there is only one unbounded type in $S_\CX(A)$, whose a
  nearest co-commutative coset is the whole $\RR^2$. Thus
  $\CL_A(\CX)=\{ \{0\}, \RR^2\}$.  However it is not hard to see that
  in any proper elementary extension $B$ of $\RR$ there are types in
  $S_\CX(B)$ whose nearest co-commutative cosets are translates of
  $L=\{0\}\times \RR$, and $\SL_B(\CX)=\{ \{0\}, L, \RR^2\}$.
\end{rem}

\begin{defn}
  For $\CX\subseteq G^\sharp$ an $\Lom$-definable set over $A$, we denote
  by $\SLm(\CX)$ the set of maximal subgroups, by inclusion, in
  $\SL_A(\CX)$. By Corollary~\ref{cor:incl-L}, it does not depend on
  $A$.
\end{defn}

We now have:
\begin{thm}\label{thm:main-section3}
  Let $G$ be an $\Rom$-definable group, $A\subseteq \mfR$, and let
  $\CX\subseteq G^\sharp$ be $\Lom$-definable over $A$.

  For every $r\in \RR^{>0}$, there are definable co-commutative
  subgroups $L_1,\ldots, L_k\subseteq G$, possibly with repetitions,
  and $a_1,\ldots,a_k\in\dcl(A)$ such that each $a_iL_i$ is a nearest
  co-commutative coset to some $p_i\in S_\CX(A)$, and
  \[ \CX \subseteq \bB_r^\sharp \ccdot\bigcup_{i=1}^k a_i\ccdot
    L_i^\sharp. \]

  In addition, every $\SLm(\CX)$ appears at least once among
  $L_1,\ldots, L_k$.
\end{thm}

\begin{proof} For each $p\in S_\CX(A)$, we choose $a_p\in \dcl(A)$
  such that $p(\mfR) \subseteq \mu\ccdot a_p\ccdot L_p^\sharp$.

  We have
  \[ \CX \subseteq \bigcup_{p\in S_\CX(A)} p(\mfR) \subseteq
    \bigcup_{p\in S_\CX(A)} \mu \ccdot a_p \ccdot L_p^\sharp \subseteq
    \bigcup_{p\in S_\CX(A)} \bB_r^\sharp\ccdot a_p \ccdot
    L_p^\sharp. \] Using logical compactness, we obtain finitely many
  $p_1,\dotsc,p_k\in S_{\CX}(A)$ such that
  \[ \CX \subseteq \bigcup_{i=1}^{k} \bB_r^\sharp\ccdot a_{p_i} \ccdot
    L_{p_i}^\sharp = \bB_r^\sharp\ccdot\bigcup_{i=1}^{k} a_{p_i}
    \ccdot L_{p_i}^\sharp . \] This proves the main part.

  In addition, let $L\in \SLm$.  Choose $p\in S_\CX(A)$ such that
  $L=L_p$ and also choose $a\in \dcl(A)$ such that
  $p(\mfR) \subseteq \mu\ccdot a \ccdot L^\sharp$.  We have
  \[ p(\mfR) \subseteq \CX \subseteq \bigcup_{i=1}^{k}
    \bB_r^\sharp\ccdot a_{p_i} \ccdot L_{p_i}^\sharp. \] Since $p$ is
  a complete type over $A$, there is $1 \leq j \leq k$ such that
  $p(\mfR) \subseteq \bB_r^\sharp\ccdot a_{p_j} \ccdot
  L_{p_j}^\sharp$.  Since $aL$ is a nearest co-commutative coset to
  $p$, by Lemma~\ref{lem:properties}, we conclude
  $L\subseteq L_{p_j}$. By maximality of $L$ we get $L=L_{p_j}$.
\end{proof}

\section{$\Gamma$-dense types in unipotent
  groups}\label{sec:gamma-dense-types}

\subsection{Preliminaries on unipotent groups}\label{sec:prel-unip-groups}

As in \cite{nilpotent}, by \emph{a unipotent group} we mean a real
algebraic subgroup of the group of real $n\times n$ upper triangular
matrices with $1$ on the diagonal.

We list below some properties of unipotent groups that we need and
refer to \cite{nilpotent} and \cite{nilpotent-book} for more details.

We fix a unipotent group $G$.

\begin{fact}\label{fact:closed}
  For a subgroup $H$ of $G$, the following are equivalent.
  \begin{enumerate}
  \item $H$ is a closed connected subgroup of $G$.
  \item $H$ is a real algebraic subgroup of $G$.
  \item $H$ is definable in $\Rom$.
  \end{enumerate}
\end{fact}

\emph{A lattice in $G$} is a discrete subgroup $\Gamma$ such that
$G/\Gamma$ is compact.

Let $\Gamma\subseteq G$ be a lattice. A real algebraic subgroup $H$ of
$G$ is called \emph{$\Gamma$-rational} if $\Gamma\cap H$ is a lattice
in $H$.

\begin{fact}\label{fact:rational-closed}
  Let $\Gamma$ be a lattice in $G$.
  \begin{enumerate}
  \item The center $Z(G)$ is $\Gamma$-rational.
  \item The commutator subgroup $[G,G]$ is closed and
    $\Gamma$-rational.
  \item If $H$ is a $\Gamma$-rational normal subgroup of $G$ and
    $\pi\colon G\to G/H$ is the quotient map then $\pi(\Gamma)$ is a
    lattice in $G/H$.  In addition, for every $\pi(\Gamma)$-rational
    subgroup $K\subseteq G/H$, the preimage $\pi^{-1}(K)$ is
    $\Gamma$-rational.
  \item If $H_1$ and $H_2$ are $\Gamma$-rational subgroups of $G$ then
    $H_1\cap H_2$ is $\Gamma$-rational as well.
  \end{enumerate}
\end{fact}

It follows from the above fact that for any real algebraic subgroup
$H$ of $G$ there is the smallest $\Gamma$-rational subgroup containing
$H$. We call it \emph{the $\Gamma$-rational closure of $H$} and denote
by $H^\Gamma$.

The next fact easily follows from Fact~\ref{fact:rational-closed}(3).
\begin{fact}\label{fact:hom and Gamma-closure} Assume that $H$ is a
  $\Gamma$-rational normal subgroup of $G$, $\pi\colon G\to G/H$ the
  quotient map and $\Gamma_0=\pi(\Gamma)$. Then for every real
  algebraic subgroup $L\subseteq G$,
  $\pi(L^\Gamma)=\pi(L)^{\Gamma_0}$.\end{fact}

We will need the following fact.
\begin{fact}\label{fact:rational-normal}
  Let $\Gamma$ be a lattice in $G$ and $H$ be a real algebraic
  subgroup of $G$.  If $H$ is a normal subgroup then $H^\Gamma$ is
  normal as well.
\end{fact}

The following is a restatement of Ratner's Orbit Closure Theorem in
the case of unipotent groups.

 \begin{fact}\cite{ratner}\label{fact:ratner}
   Let $\Gamma$ be a lattice in $G$ and $H$ be a real algebraic
   subgroup of $G$. The topological closure of $H\ccdot \Gamma$ in $G$
   is $H^\Gamma\ccdot \Gamma$.
 \end{fact}

 We will be using the following well-known fact.

\begin{fact}\label{fact:commens}
  Let $H$ be a real algebraic subgroup of $G$ and $\Gamma_1,\Gamma_2$
  be lattices in $G$. If $\Gamma_1$ and $\Gamma_2$ are commensurable,
  i.e. $\Gamma_1\cap \Gamma_2$ is of finite index in both $\Gamma_1$
  and $\Gamma_2$, then $H^{\Gamma_1}=H^{\Gamma_2}$.
\end{fact}

\subsection{$\Gamma$-dense sets in unipotent groups.}\label{sec:gamma-dense-subsets}

Let $G$ be a unipotent group and $\Gamma\subseteq G$ be a lattice.  We
say that a subset $X \subseteq G$ is \emph{$\Gamma$-dense in $G$} if
the set $X\ccdot\Gamma$ is dense in $G$, i.e. $\cl(X\ccdot\Gamma)=G$.
Using Lemma~\ref{lem:st-closed}(2), we conclude that a subset
$X\subseteq G$ is $\Gamma$-dense in $G$ if and only if
$\st(X^\sharp\ccdot \Gamma^\sharp)=G$. We use this fact to extend the
notion of $\Gamma$-density to arbitrary subsets of $G^\sharp$.

\begin{defn} Let $G$ be a unipotent group, $\Gamma\subseteq G$ a
  lattice and $\CX\subseteq G^\sharp$ be an arbitrary set.
  \begin{enumerate}
  \item We say that $\CX$ is \emph{$\Gamma$-dense in $G$} if
    $\st(\CX\ccdot \Gamma^\sharp)= G$.
  \item We say that $\CX$ is \emph{strongly $\Gamma$-dense in $G$} if
    $\st(\CX\ccdot \Gamma_1^\sharp)= G$ for every lattice $\Gamma_1$
    commensurable with $\Gamma$.
  \item We say that a type $p\in S_G(A)$ is \emph{(strongly)
      $\Gamma$-dense in $G$} if the set $p(\mfR)$ is (strongly)
    $\Gamma$-dense in $G$.
  \end{enumerate}
\end{defn}
\begin{rem} Let $G$ is be a unipotent group, $\Gamma \subseteq G$ a
  lattice and $\CX\subseteq G^\sharp$. It is easy to see that $\CX$ is
  strongly $\Gamma$-dense in $G$ if and only if it is $\Gamma_0$-dense
  for every subgroup $\Gamma_0\subseteq\Gamma$ of finite index.
\end{rem}

\begin{sample}
  Let $G=(\RR,+)$, $\Gamma=\ZZ$ and let $X$ be the closed interval
  $[0,1]$. The set $X^\sharp$ is $\Gamma$-dense in $G$, but not
  strongly $\Gamma$-dense.
\end{sample}

The following fact follows from Facts~\ref{fact:ratner} and
\ref{fact:commens}.

\begin{fact}\label{factdense-subgrop}
  Let $G$ be a unipotent group and $L\subseteq G$ a real algebraic
  subgroup. For a lattice $\Gamma\subseteq G$ the following are
  equivalent.
  \begin{enumerate}
  \item $L$ is $\Gamma$-dense in $G$.
  \item $L^\Gamma=G$
  \item $L$ is strongly $\Gamma$-dense in $G$.
  \end{enumerate}

\end{fact}

We observe:
\begin{lem}\label{lem:Gamma-dense} Let $\Gamma$ be a lattice in a
  unipotent group $G$. A subset $\CX\subseteq G^\sharp$ is $\Gamma$-dense
  in $G$ if and only if $\mu\cdot \CX\cdot \Gamma^\sharp=G^\sharp$.
\end{lem}
\begin{proof} The ``if'' part is clear.

  For ``the only if'' part, since $G/\Gamma$ is compact, given
  $g\in G^\sharp$ there is $\gamma\in \Gamma^\sharp$ such that
  $g\gamma\in \CO$. Thus, since $\CX$ is $\Gamma$-dense in $G$, there
  is $a\in \CX$ such that $\st(g\gamma)= \st(a)$.  It follows that
  $g\in \mu \ccdot a \ccdot \Gamma^\sharp$.
\end{proof}

We will need the following fact.

\begin{fact}[\cite{nilpotent}*{Lemma~5.1}]\label{fact:51}
  Let $\pi\colon G\to H$ be a real algebraic surjective homomorphism
  of unipotent groups, and $\CX$ a subset of $G^\sharp$.  Then, for
  every lattice $\Gamma \subseteq G$, if $\pi(\Gamma)$ is closed in
  $H$ then
  \[\pi(\st(\CX\ccdot \Gamma^\sharp)) = \st(\pi (\CX )\ccdot
    \pi(\Gamma^\sharp)).\]
\end{fact}
Our main goal is to describe $\Gamma$-dense types. We will consider
the abelian case first.

\subsection{$\Gamma$-dense types in abelian groups.}\label{sec:gamma-dense-types-1}

Since every abelian unipotent group is algebraically isomorphic to
$(\RR^m,+)$ for some $m$, we often identify an abelian unipotent group
with an $\RR$-vector space $(\RR^m,+)$.

In the abelian case, every subgroup is co-commutative, hence for a set
$A\subseteq \mfR$ and a type $p(x)\in S(A)$ on $\RR^m$, instead of a
nearest co-commutative coset to $p$ we say \emph{a nearest coset to
  $p$}.

Our first goal of this section is to prove the following:
\begin{prop}\label{prop:dense types} Let $G$ be an abelian unipotent
  group, $A\subseteq \mfR$, $p\in S_G(A)$, and $a_p\in \dcl(A)$ be
  such that $a_p+L_p$ is a nearest coset to $p$. For a lattice
  $\Gamma\subseteq G$, the following are equivalent.
  \begin{enumerate}
  \item The type $p$ is $\Gamma$-dense in $G$.
  \item $L_p^\Gamma=G$.
  \item The type $p$ is strongly $\Gamma$-dense in $G$.
  \end{enumerate}
\end{prop}

\begin{proof}

  Since, by Fact~\ref{fact:commens}, $L^\Gamma=L^{\Gamma_0}$ for any
  subgroup $\Gamma_0\subseteq \Gamma$ of finite index, it is
  sufficient to show $(1) \Leftrightarrow (2)$.

  For simplicity we denote $L_p^\Gamma$ by $L$, and assume
  $A=\dcl(A)$.

  \medskip
  \noindent$(1)\Rightarrow (2)$.  Assume that $L \neq G$, hence, by
  Fact~\ref{fact:ratner}, the set $L+\Gamma$ is a closed proper
  subgroup of $G$. By Lemma~\ref{lem:st-closed}(2),
  $\mu+L^\sharp+\Gamma^\sharp$ is a proper subgroup of $G^\sharp$,
  hence the coset $a_p +\mu+ L^\sharp +\Gamma^\sharp$ is a proper
  subset of $G^\sharp$.

  Since
  $\mu+p(\mfR)+\Gamma^\sharp \subseteq a_p +\mu+ L^\sharp
  +\Gamma^\sharp$, the set $\mu+p(\mfR)+\Gamma^\sharp$ is a proper
  subset of $G^\sharp$ and, by Fact~\ref{lem:Gamma-dense}, $p(\mfR)$
  is not $\Gamma$-dense, so (1) fails.

  \medskip

  \noindent$(2)\Rightarrow (1)$. Assume $L=G$ and we prove that
  $\st(p(\mfR)+ \Gamma^\sharp)=G$, The proof is similar to
  \cite{nilpotent}*{Proposition 5.3}.

  We use induction on $\dim G$.

  If $\dim(G)=0$ then there is nothing to prove.

  Assume $\dim(G)=n > 0$ and the result holds for all abelian
  unipotent groups of dimension less than $n$.

  We have $\dim(L)>0$, hence, by Remark~\ref{rem:nearshort}, $p$ is a
  long type.  Let $P$ be the $\mu$-stabilizer of $p$.  By
  Proposition~\ref{prop:main-mustab}, $P$ is a real algebraic subgroup
  of $G$ of positive dimension, and by Lemma~\ref{lem:mu-nearest},
  $P\subseteq L_p$, hence $P^\Gamma\subseteq L_p^\Gamma=L$.

  Let $\pi:G\to G_0 := G/P^\Gamma$ be the quotient map,
  $\Gamma_0=\pi(\Gamma)$, and $q=\pi(p)$.  By
  Fact~\ref{fact:rational-closed}, $\Gamma_0$ is a lattice in
  $G_0$. Notice that $\dim(G_0)<\dim(G)$.

  Let $a_q+L_q$ be a nearest coset to $q$. It follows from
  Lemma~\ref{lem:near-hom} that $L_q=\pi(L_p)$. Since $G=L_p^\Gamma$,
  by Fact~\ref{fact:rational-closed}(3), $G_0=L_q^{\Gamma_0}$, hence,
  by induction hypothesis, the type $q$ is $\Gamma_0$-dense in
  $G_0^\sharp$, and
  \[\st(q(\mfR) +\Gamma_0^\sharp) = G_0. \]
  Applying Fact~\ref{fact:51}, we obtain
  \[
    \pi(\st(p(\mfR)+\Gamma^\sharp))=\st(q(\mfR)+\Gamma_0^\sharp)=G_0. \]

  Let $D=\st(q(\mfR) +\Gamma^\sharp)$.  By
  Lemma~\ref{lem:st-closed}(3), it is a closed subset of $G$.  It is
  not hard to see that $D$ is invariant ander the action of both $P$
  and $\Gamma$, hence it is invarint under $P+\Gamma$. Since $D$ is
  closed, it is invariant under the action of the topological closure
  of $P+\Gamma$.  By Fact~\ref{fact:ratner}, $P^\Gamma$ is contained
  in $\cl(P+\Gamma)$, hence $D$ is invariant under $P^\Gamma$.  Since
  $\pi(D)=G_0$ and $\ker(\pi)=P^\Gamma$, it follows then $D=G$, hence
  $p$ is $\Gamma$-dense in $G$.

  This finishes the proof of Proposition~\ref{prop:dense types}.
\end{proof}

As a corollary we obtain the following theorem.

\begin{thm}\label{thm:main types} Let $A\subseteq \mfR$,
  $G=(\RR^n,+)$, $p\in S_G(A)$, and let $a_p+L_p$ be a nearest coset
  to $p$ (so $a_p\in \dcl(A)$).  Then for every lattice
  $\Gamma \subseteq \RR^n$, we have
  \[ \mu+p(\mfR)+\Gamma^\sharp= \mu+a_p +L_p^\sharp +\Gamma^\sharp .\]
\end{thm}
\begin{proof}

  Consider the type $p_1=-a_p+p$. It is a complete $\Lom$-type over
  $A$.  Clearly $L_p$ is a nearest coset to $p_1$, hence there exists
  a type $p_2\in S_{L_p}(A)$ which is $\mu$-equivalent to $p_1$, and
  therefore $L_p$ is also a nearest coset to $p_2$.  Let
  $G_0=L_p^\Gamma$ and $\Gamma_0=\Gamma\cap G_0$, a lattice in $G_0$.

  Working in $G_0$ we have that $L_p^{\Gamma_0}=G_0$, hence by
  Proposition~\ref{prop:dense types}, the type $p_2$ is
  $\Gamma_0$-dense in $G_0$, so, by Lemma~\ref{lem:Gamma-dense},
  \[(\mu\cap G_0^\sharp)+p_2(\CR)+\Gamma_0^\sharp=G_0^\sharp.\]

  Obviously, $L_p$ is also $\Gamma_0$-dense in $G_0$, hence
  $G_0^\sharp=(\mu\cap G_0^\sharp)+L_p^\sharp+\Gamma_0^\sharp$.

  We conclude
\[\mu+p(\CR)+\Gamma^\sharp=\mu+a_p+\mu+p_2(\CR)+\Gamma^\sharp=\mu+a_p+L_p^\sharp+\Gamma^\sharp.\]

\end{proof}

\subsection{Abelianization and density}\label{sec:abelinization}
For a unipotent group $G$ we will denote by $\Gab$ the abelianization
of $G$, i.e. the group $\Gab= G/[G,G]$, and by $\piab$ the quotient
map $\piab \colon G\to \Gab$. The group $\Gab$ is also unipotent and
$\dim G>0$ if and only if $\dim\Gab >0$.

If $\Gamma \subseteq G$ is a lattice then we denote by $\Gamab$ the
group $\Gamab=\piab(\Gamma)$. By Fact~\ref{fact:rational-closed},
$\Gamab$ is a lattice in $\Gab$. Our main goal in this section to show
that a type $p\in S_G(A)$ is $\Gamma$-dense in $G$ if and only if its
abelianization $\piab(p)$ is $\Gamab$-dense in $\Gab$.

The next proposition is a key.

  \begin{prop}\label{prop:ifpart}
    Let $G$ be a unipotent group, $A\subseteq \mfR$, $\Gamma$ a
    lattice in $G$ and $p\in S_G(A)$. Assume $p$ is not $\Gamma$-dense
    in $G$. Then there is a co-commutative $\Gamma$-rational subgroup
    $H\subseteq G$ such that for the projection $\pi\colon G\to G/H$
    the type $\pi(p)$ is not $\pi(\Gamma)$-dense in $G/H$.
  \end{prop}
  \begin{proof}
    By induction on $\dim(G)$.

    If $\dim(G)=0$ then there is nothing to prove.

    Assume $\dim(G)=n > 0$ and the proposition holds for all unipotent
    groups of dimension less than $n$.

    If the type $p$ is short then, by Lemma~\ref{lem:short1}(1), the
    type $\piab(p)$ is short as well, and it is easy too see,
    e.g. using Lemma~\ref{lem:short2}(1), that a short type is not
    $\Gamab$-dense in $\Gab$. We can take $H=[G,G]$ that is
    $\Gamma$-rational by Fact~\ref{fact:rational-closed}(2).

    Thus we may assume that $p$ is a long type. Let $P$ be the
    $\mu$-stabilizer of $p$.  By Proposition~\ref{prop:main-mustab},
    $P$ is an $\Rom$-definable subgroup of $G$ of positive dimension.

    As in \cite[Propositioni 5.3]{nilpotent}, we consider the smallest
    $\Rom$-definable, normal $\Gamma$-rational subgroup of $G$
    containing $P$ and denote it by $N(P)^\Gamma$.  Let $N_0$ be the
    intersection of $N(P)^\Gamma$ with the center of $G$. Since $G$ is
    unipotent and $N(P)^\Gamma$ has positive dimension, the group
    $N_0$ is also of positive dimension (see, for example,
    \cite[Proposition 7.13]{stroppel}), and, by
    Fact~\ref{fact:rational-closed}, it is $\Gamma$-rational.

    Let $\pi\colon G\to G_0 := G/N_0$ be the quotient map,
    $\Gamma_0=\pi_0(\Gamma)$, and $q=\pi(p)$.  By
    Fact~\ref{fact:rational-closed}, $\Gamma_0$ is a lattice in $G_0$.

    \medskip

    We claim that the type $q$ is not $\Gamma_0$-dense in $G_0$.

    \medskip

    Indeed, assume towards contradiction that $q$ is $\Gamma_0$-dense
    in $G_0$. Then, by Fact~\ref{fact:51},
    \[\pi_0( \st(p(\mfR) \ccdot \Gamma^\sharp))=G_0.\]
    Let $D_{p,\Gamma}=\st(p(\mfR)\ccdot\Gamma^\sharp)$ It follows from
    the above equation that
    \begin{equation}
      \label{eq:4}
      D_{p,\Gamma}\ccdot N_0 =G.
    \end{equation}

    Our aim is to show that $D_{p,\Gamma} =G$, contradicting the fact
    that $p$ is not $\Gamma$-dense in $G$.

    Since, by Lemma~\ref{lem:st-closed}(2), the set $D_{p,\Gamma}$ is
    a closed subset of $G$, it is sufficient to show that it is dense
    in $G$.

    \medskip
    \noindent\textbf{Claim A.} {\itshape The set $D_{p,\Gamma}$ is
      left invariant under the $\mu$-stabilizer $P$ of $p$.}
    \begin{proof}
      Note that
      $D_{p,\Gamma}=\st(\mu\ccdot p(\mfR)\ccdot\Gamma^\sharp)$, and
      $\mu p$ is left-invariant under $P$. Thus, for $g\in P$,
  \[g\cdot D_{p,\Gamma}=g\ccdot\st(\mu\ccdot p(\mfR)\ccdot\Gamma^\sharp)=\st(g\ccdot\mu
  \ccdot p(\mfR)\ccdot \Gamma^\sharp)=\st(\mu \ccdot p(\mfR)\ccdot
  \Gamma^\sharp)=D_{p,\Gamma}.\]
\end{proof}

Clearly $D_{p,\Gamma}$ is right-invariant under action of $\Gamma$.
Thus, $P\ccdot D_{p,\Gamma}\ccdot \Gamma= D_{p,\Gamma}$, and, in
addition, by equation~\eqref{eq:4}, $D_{p,\Gamma}N_0=G$.

Because $N_0=N(P)^\Gamma\cap Z(G)$, our goal, $D_{p,\Gamma}=G$,
follows from the following general result:

\medskip
\noindent\textbf{Claim B.} {\itshape For a unipotent group $G$, assume
  that $D\subseteq G$ is a closed set, left invariant under a real
  algebraic subgroup $P\subseteq G$ and right invariant under a
  lattice $\Gamma\subseteq G$.  Let $N_0\subseteq N(P)^\Gamma\cap N_G(P)$ be a
  subgroup of $G$. If $DN_0=G$ then $D=G$.}
\begin{proof} Let $Y=\{g\in G: (P^g)^\Gamma=N(P)^\Gamma\}$. This is
  not, in general, a definable set, but, by
  \cite{nilpotent}*{Proposition~4.3}, it is dense in $G$.  Thus, it is
  sufficient to prove that $Y\subseteq D$.

  First note that $Y$ is left invariant under $N_0$. Indeed, assume
  that $a\in N_0b$. Since $ab^{-1}\in N_0\subseteq N_G(P)$, then $P^a=P^b$,
  implying that $b\in Y$ if and only if $a\in Y$.

  Let $b\in Y$. Since $DN_0=G$, there is $a\in D$ such that
  $b\in aN_0$, and therefore $a\in Y$.  Using the definition of $Y$
  and the fact that $aP^a=Pa$, we obtain
\[b\in a\ccdot N_0\subseteq a\ccdot  N(P)^\Gamma=a\ccdot (P^a)^\Gamma=a\cdot
\cl(P^a\ccdot \Gamma)\subseteq \cl(a \ccdot P^a \ccdot
\Gamma)=\cl(P\ccdot a \ccdot \Gamma).\]
Since $a\in D$, by the
invariance properties of $D$, also $b\in D$. Hence $Y\subseteq D$, so
$D=G$, a contradiction.

This ends the proof of Claim B, and thus we conclude that $q$ is not
$\Gamma_0$-dense in $G_0$.  \end{proof}

Applying the induction hypothesis to $G_0$, $\Gamma_0$ and $q$, we
obtain a co-commutative $\Gamma_0$-rational subgroup
$H_0\subseteq G_0$ such that the image of $q$ in $G_0/H_0$ is not
$\Gamma_0/H_0$-dense.  It is not hard to see that $H=\pi_0^{-1}(H_0)$
is a co-commutative $\Gamma$-rational subgroup of $G$ satisfying the
conclusion of the proposition.

This finishes the proof of Proposition~\ref{prop:ifpart}

\end{proof}

We can now prove the main theorem of this section.
\begin{thm}\label{thm:abel}
  Let $G$ be a unipotent group, $A\subseteq \mfR$ and $\Gamma$ a lattice in
  $G$.  A type $p\in S_G(A)$ is $\Gamma$-dense in $G$ if and only if
  the type $\piab(p)$ is $\Gamab$-dense in $\Gab$.
\end{thm}
\begin{proof} Let $q=\piab(p)$. We write additively the group
  operation in $\Gab$. By Fact~\ref{fact:51},
  \[ \piab(\st(p(\mfR)\cdot \Gamma^\sharp)) =
    \st(q(\mfR)+\Gamab^\sharp).\] This implies the ``only if'' part.

  We prove the ``if part'' part by contraposition, using
  Proposition~\ref{prop:ifpart}.  Indeed, assume the type $p$ is not
  $\Gamma$-dense in $G$, and we derive that $\piab(p)$ is not
  $\Gab$-dense in $\Gab$.

  Let $H$ and $\pi\colon G\to G/H$ be as in
  Proposition~\ref{prop:ifpart}.  Since $H$ contains $[G,G]$, the map
  $\pi$ factors through $\Gab$, i.e. there is $\pi'\colon \Gab\to G/H$
  with $\pi=\pi'\co \piab$.  By Fact~\ref{fact:51}
  \[\pi'\Bigl(\st\bigl(\piab(p)(\mfR)+\Gamab^\sharp\bigr)\Bigr)
    = \st\Bigl( \pi(p)(\mfR)+\pi(\Gamma)^\sharp\Bigr).\] Since
  $\pi(p)$ is not $\pi(\Gamma)$-dense in $G/H$, the type $\piab(p)$ is
  not $\Gamab$-dense in $\Gab$.
\end{proof}

The following summarizes main results of this section.

\begin{thm}\label{thm:type-dense} Let $G$ be a unipotent group,
  $A\subseteq \mfR$ and $p\in S_G(A)$.

  For a lattice $\Gamma\subseteq G$, the following are equivalent.

  \begin{enumerate}

  \item The type $p$ is $\Gamma$-dense in $G$.
  \item The type $p$ is strongly $\Gamma$-dense in $G$.
  \item $L_p^\Gamma=G$.
  \item The type $\piab(p)$ is $\Gamab$-dense in $\Gab$.
  \item The type $\piab(p)$ is strongly $\Gamab$-dense in $\Gab$.
  \item $\piab(L_p)^{\Gamab}=\Gab$.
  \end{enumerate}
\end{thm}

\section{$\Gamma$-dense definable subsets of unipotent groups}\label{sec:gamma-dense-defin}
We fix a unipotent group $G$.

In this section we obtain a description of $\Gamma$-dense definable
subsets of $G$ similar to that of Theorem~\ref{thm:type-dense}.

First an elementary lemma.

\begin{lem}\label{lem:elem1}
  For $A\subseteq \mfR$, let $\CX\subseteq G^\sharp$ be
  $\Lom$-definable over $A$. For a lattice $\Gamma\subseteq G$, if
  some type $p\in S_\CX(A)$ is $\Gamma$-dense in $G$ then $\CX$ is
  strongly $\Gamma$-dense $G$.
\end{lem}
\begin{proof} Assume a type $p\in S_\CX(A)$ is $\Gamma$-dense, hence,
  by Theorem~\ref{thm:type-dense}, it is strongly
  $\Gamma$-dense. Since $p(\mfR)\subseteq \CX$, obviously $\CX$ is
  strongly $\Gamma$-dense.
\end{proof}

The main goal of this section (Theorem~\ref{thm:main1}) is to show
that an appropriate converse of the above lemma holds. Namely, an
$\Lom$-definable subset $\CX\subseteq G^\sharp$ is strongly
$\Gamma$-dense in $G$ if and only if some type on $\CX$ is
$\Gamma$-dense.

We need two lemmas and a proposition.

 \begin{lem}\label{lem:fst-cyl}
   Let $L\subseteq G$ be a normal real algebraic subgroup and
   $\alpha\in G^\sharp$.
   \begin{enumerate}
   \item For any lattice $\Gamma\subseteq G$ the set
     $\CO\cap \alpha\Gamma^\sharp$ is nonempty, and for every
     $g\in \st( \alpha\Gamma^\sharp)$ we have
     \[ \st(\alpha\ccdot L^\sharp\ccdot \Gamma^\sharp)= g \ccdot
       L^\Gamma\ccdot \Gamma. \]
   \item If $L$ is co-commutative then $\alpha L$ is a nearest
     co-commutative coset to some type $p\in S(\alpha)$ on
     $\alpha L^\sharp$, hence $\SLm(\alpha L^\sharp)=\{ L\}$.
   \end{enumerate}
 \end{lem}
 \begin{proof}
   (1) Because $\Gamma$ is co-compact, the set
   $\CO\cap \alpha\Gamma^\sharp$ is not empty. Let
   $g\in \st( \alpha\Gamma^\sharp)$, and we choose
   $\gamma\in\Gamma^\sharp$ with $\alpha \cdot\gamma\in \mu g$.  Since
   $L$ is normal, we have
   \begin{multline*}
     \st(\alpha\ccdot L^\sharp \ccdot \Gamma^\sharp)
     =\st(\alpha\ccdot\Gamma^\sharp\ccdot L^\sharp) =\st(g\ccdot
     \Gamma^\sharp\ccdot L^\sharp)
     \\
     =\st(g\ccdot L^\sharp\ccdot \Gamma^\sharp) =\cl(g\ccdot L \ccdot
     \Gamma) =g\ccdot L^\Gamma\ccdot \Gamma.
   \end{multline*}

   (2) We first consider the abelian case, namely we assume
   $G=(\RR^n ,+)$ and $L\subseteq \RR^n$ a linear subspace.  We need to
   show that there is $\beta\in \alpha +L^\sharp$ such that
   $\beta\notin \mu (\gamma +L_0^\sharp)$, for any
   $\gamma\in \dcl(\alpha)$ and a proper subspace $L_0 \subseteq L$.
   It is thus sufficient to show
   \begin{equation} \label{eq:1} \alpha+L^\sharp \not\subseteq \bigcup
     \{\bB_1^\sharp+\gamma+L_0^\sharp \colon \gamma\in \dcl(\alpha),
     L_0\subsetneq L \text{ a subspace}\}.
   \end{equation}

   By logical compactnes,s \eqref{eq:1} follows from the following
   claim.  \medskip

\noindent{\bf Claim.} {\itshape Let $r\in \RR^{\geq 0}$, and
  $L_1,\dotsc,L_k\subseteq L$ be proper subspaces. Then for any
  $\alpha,\gamma_1,\dotsc,\gamma_k\in \mfR^n$ we have
  \[ \alpha+L^\sharp \not\subseteq \bB_r^\sharp+\bigcup_{i=1}^k
    \bigl(\gamma_i+L_i^\sharp\bigr). \]}

\begin{proof}[Proof of Claim.]
  Since, for fixed $r$ and $L_1,\dotsc,L_k$, the conclusion of the
  claim can be expressed by a first-order formula, we can work in
  $\RR$ instead of $\mfR$;and also, subtracting $\alpha$ from both
  sides, we only need to consider the case $\alpha=0$.

  We fix proper subspaces $L_1,\dotsc,L_k\subseteq L$, and show that
  for all $\geq 0$ and $b_1,\dotsc,b_k\in \RR^n$, we have
  $L\not\subseteq \bB_r+\bigcup_{i=1}^k(b_i+L_i)$.

  Clearly, by the dimension assumptions,
  $L \neq L_1\cup\dotsc\cup L_k$.

  Next, let us see that $L\not\subseteq \bB_r +\bigcup_{i=1}^k L_i$
  for any $r\in \RR^{\geq 0}$. Indeed, choose
  $c\in L \setminus (\bigcup_{i=1}^k L_i)$.  Then, for every
  $i=1,\dotsc,k$, we have $d(c,L_i)>0$, where $d(\cdot,\cdot)$ denotes
  the Euclidean distance in $\RR^n$.  Since $d(tc,L_i)=t d(c,L_i)$ for
  $t\in \RR^{>0}$, we obtain that for given $r\in \RR^{\geq}$ and
  $i=1,\dotsc,k$, for $t$ large enough, $d(tc,L_i)>r$, and hence
  $tc \in L\setminus (\bB_r +\bigcup_{i=1}^k L_i)$.

  Finally, assume that for some $r>0$ and $b_1,\dotsc,b_k\in \RR^n$ we
  would have $L\subseteq \bB_r+\bigcup_{i=1}^k(b_i+L_i)$. Then,
  choosing $r'>0$ big enough so that $b_i+\bB_r \subseteq \bB_{r'}$
  for $i=1,\dotsc,k$, we would have
  $L\subseteq \bB_{r'} +\bigcup_{i=1}^k L_i$, a contradiction.

  This finishes the proof of Claim, and hence the lemma, in the case
  that $G$ is abelian.

  When $G$ is nilpotent and $L\subseteq G$ is co-commutative we first apply
  the above to $\piab(\alpha L^\sharp)\subseteq \Gab^\sharp$ and find a
  type over $\alpha$ with $q\vdash \piab(\alpha L^\sharp)$, such that
  $\piab(\alpha L)$ is a nearest coset to $q$. Now, choose a type $p$
  over $\alpha$ such that $p\vdash \alpha L$ and $\piab(p)=q$. A
  nearest co-commutative coset to $p$ is contained in $\alpha L$, and
  projects via $\piab$ onto $\piab(\alpha L)$ (see Lemma
  \ref{lem:near-hom}). Since $L\supseteq [G,G]$, it follows that
  $\alpha L$ is a nearest co-commutative coset to $p$. This finishes
  the proof of the claim.
\end{proof}
End of the proof of the lemma.

\end{proof}

We are going to need the following result.

\begin{lem}\label{lem:proper}
  If $H$ is a proper real algebraic subgroup of $G$ then $\piab(H)$ is
  a proper subgroup of $\Gab$
\end{lem}
\begin{proof}
  By \cite{nilpotent-book}*{Theorem~1.1.13} there is a chain of real
  algebraic subgroups
  \[ \{e\}=H_0\subsetneq \dotsb \subsetneq H=H_m \subsetneq
    H_{m+1}\subsetneq \dotsb \subsetneq H_n=G, \] with $n=\dim(G)$ and
  $\dim H_{i+1}=\dim H_i+1$.  By \cite{nilpotent-book}*{Lemma~1.1.8},
  $[G,G]\subseteq H_{n-1}$.  Hence $H_{n-1}/[G,G]$ is a proper
  subgroup of $\Gab$ and so is $H/[G,G]$.
\end{proof}

\begin{prop}\label{prop:finite-index}
  Let $\Gamma\subseteq G$ be a lattice, and $L_1,\ldots, L_k$ proper
  $\Gamma$-rational subgroups of $G$. If $K\subset G$ is a compact set
  then there is a subgroup $\Gamma_0\subseteq \Gamma$ of finite index such
  that for any $g_1,\ldots, g_k\in G$, we have
  \[
    K\ccdot \bigcup_{i=1}^k g_i\ccdot L_i \ccdot \Gamma_0\neq G.\]
\end{prop}

\begin{proof}
  We first consider the case when $G$ is abelian.  So we assume
  $G=(\RR^n,+)$.

  \medskip
  \noindent\textbf{Claim.}
  \emph{Let $\Gamma \subseteq \RR^n$ be a lattice, $K\subseteq \RR^n$
    a compact set, and $L\subseteq \RR^n$ be a proper
    $\Gamma$-rational subspace.  Then for any $m\in \NN$, there is a
    subgroup $\Gamma'\subseteq \Gamma$ of finite index such that for
    some $b_1,\dotsc,b_m \in \RR^n$, the translates $b_i+K+L+\Gamma'$,
    $i=1,\dotsc,m$, are pairwise disjoint.}
  \begin{proof}
    Replacing $\RR^n$ by $\RR^n/L$ if needed, we may assume that $L$
    is the trivial subspace $\{0\}$.

    Since $K$ is compact, it is bounded, hence there are
    $b_1,\dotsc b_m\in \RR^n$ such that the translates
    $b_1+K, \dotsc, b_m+K$ are pair-wise disjoint.  Let
    $B=\bigcup_{i=1}^m (b_i+K)$. Obviously $B$ is compact and hence
    the set $B'=B-B=\{b-b'\colon b,b'\in B\}$ is compact as well.

    Since $\Gamma$ is discrete, the intersection $\Gamma\cap B'$ is
    finite. Every finitely generated abelian group is residually
    finite, i.e. the intersection of all subgroups of finite index is
    trivial, hence there is a subgroup $\Gamma'\subseteq \Gamma$ of
    finite index with $\Gamma'\cap B'=\{0\}$.  It is not hard to see
    that the sets $b_i+K+\Gamma'$, $i=1,\dotsc m$, are pairwise
    disjoint.

    This finishes the proof of the claim. \end{proof}

  We return to the proof of the proposition for $G=(\RR^n,+)$.  We
  apply the above claim to each $L_i$ with $m=k+1$, and for each
  $i=1,\dotsc,k$, obtain a subgroup $\Gamma_i\subseteq G$ of finite
  index such that $K+L_i+\Gamma_i$ has $k+1$ disjoint translates.

  Since every abelian group is amenable, there is a $G$-invariant
  finitely additive probability measure
  $\lambda\colon \mathscr{P}(G)\to [0,1]$.  By our choice of
  $\Gamma_i$, we have $\lambda(K+L_i+\Gamma_i)\leq 1/(k+1)$.

  We take $\Gamma_0=\bigcap_{i=1}^k \Gamma_i$. For any
  $g_1,\dotsc,g_k\in G$ we have
  \begin{multline*}
    \lambda \Bigl(\bigcup_{i=1}^k g_i+K+ L_i+\Gamma_0\Bigr )\leq
    \sum_{i=1}^k \lambda( g_i+K+ L_i+\Gamma_0) \leq\\
    \sum_{i=1}^k \lambda( g_i+K+ L_i+\Gamma_i)\leq k/(k+1)<1.
  \end{multline*}
  Hence $\bigcup_{i=1}^k g_i+K+ L_i+\Gamma_0\neq G$.  This finishes
  the proof of the abelian case.

  \medskip Assume now that $G$ is an arbitrary unipotent group. Let
  $K_\mathrm{ab}=\piab(K)$, $\Gamab=\piab(\Gamma)$, and, for
  $i=1,\dotsc,k$, let $L_i^\mathrm{ab}=\piab(L_i)$.  Obviously
  $K_\mathrm{ab}$ is a compact subset of $\Gab$, $\Gamab$ is a lattice
  in $\Gab$ by Fact~\ref{fact:rational-closed}, and it is not hard to
  see that each $L_i^\mathrm{ab}$ is $\widetilde{\Gamma}$-rational
  subgroup of $\Gab$.  It also follows from Lemma~\ref{lem:proper}
  that each $L_i^\mathrm{ab}$ is a proper subgroup of $\Gab$.

  We now use the abelian case and find a subgroup
  $\Gamma_0'\subseteq \Gamab$ such that for any
  $b_1,\dotsc,b_m\in \Gab$ we have
  $ K_\mathrm{ab}\ccdot\bigcup_{i=1}^k b_i\ccdot L_i^\mathrm{ab}\ccdot
  \Gamma_0'\neq \Gab$.

  We take $\Gamma_0=\piab^{-1}(\Gamma_0') \cap \Gamma$.
\end{proof}

We are now ready to prove one of the main theorems of this paper.

\begin{thm}\label{thm:main1} Let $G$ be a unipotent group,
  $A\subseteq \mfR$, and let $\CX\subseteq G^\sharp$ be a set
  $\Lom$-definable over $A$.

  For a lattice $\Gamma\subseteq G$, the following are equivalent:
  \begin{enumerate}[(a)]
  \item The set $\CX$ is strongly $\Gamma$-dense in $G$.
  \item $L^\Gamma=G$ for some $L\in\SLm(\CX)$.
  \item Some type $p\in S_\CX(A)$ is $\Gamma$-dense.
  \end{enumerate}
\end{thm}

\begin{proof}
  By Theorem~\ref{thm:type-dense}, $(b)\Leftrightarrow (c)$, and, by
  Lemma~\ref{lem:elem1}, $(c)\Rightarrow (a)$.

  \medskip

  Let us show that $(a)\Rightarrow (b)$.

  Let $\Gamma\subseteq G$ be a lattice.  We choose $L_i$ and $a_i$,
  $i=1, \dotsc, k$, as in Theorem~\ref{thm:main-section3} with $r=1$.

  Assume $(b)$ fails, namely, $L^\Gamma\neq G$ for all
  $L\in \SLm(\CX)$. Then clearly, $L_i^\Gamma\neq G$, for all
  $i=1,\dotsc,k$. For any subgroup $\Gamma_0 \subseteq G$ of finite
  index we have
  \[\CX \ccdot \Gamma_0^\sharp
    \subseteq \bB_1^\sharp \ccdot \bigcup_{i=1}^k a_i \ccdot
    L_i^\sharp\ccdot \Gamma_0^\sharp ,\] hence
  \[\st(\CX \ccdot \Gamma_0^\sharp)
    \subseteq \st\Bigl( \bB_1^\sharp\ccdot \bigcup_{i=1}^k a_i\ccdot
    L_i^\sharp \ccdot \Gamma_0^\sharp\Bigr) =\bigcup_{i=1}^k
    \st(\bB_1^\sharp\ccdot a_i\ccdot L_i^\sharp \ccdot
    \Gamma_0^\sharp).\]

  Using Lemma~\ref{lem:fst-cyl}(1) we choose $g_1,\dotsc,g_k\in G$
  such that
  $\st(\bB_r^\sharp\ccdot a_i\ccdot L_i^\sharp \ccdot \Gamma_0^\sharp)
  =\bB_r\ccdot g_i \ccdot L_i^{\Gamma_0}\ccdot \Gamma_0$. By
  fact~\ref{fact:commens}, $L_i^{\Gamma_0}=L_i^\Gamma$, hence
  \[ \st(\CX \ccdot \Gamma_0^\sharp)\subseteq \bB_r\ccdot
    \bigcup_{i=1}^k g\ccdot L_i^\Gamma\ccdot \Gamma_0.\]

  By Proposition~\ref{prop:finite-index}, there exists
  $\Gamma_0\subseteq \Gamma$ of finite index for which the set on the
  right is a proper subset of $G$, hence $(a)$ fails. Thus,
  $(a)\Rightarrow (b)$.

\end{proof}

Recall that for $\piab:G\to \Gab$, and $\Gamma\subseteq G$ a lattice, we
let $\Gamab=\piab(\Gamma)$.
\begin{cor}\label{cor:xiffxab}
  Let $G$ be a unipotent group, and let $\CX\subseteq G^\sharp$ be an
  $\Lom$-definable set.  For a lattice $\Gamma\subseteq G$, the set
  $\CX$ is strongly $\Gamma$-dense in $G$ if and only if $\piab(\CX)$
  is strongly $\Gamab$-dense in $\Gab$.
\end{cor}
\begin{proof}
  The ``only if'' part follows from Fact~\ref{fact:51}.

  For the ``if'' part, assume that $\piab(\CX)$ is strongly
  $\Gamab$-dense in $\Gab$. Choose a set $A\subseteq \mfR$ such that
  $\CX$ is $\Lom$-definable over $A$. Applying
  Theorem~\ref{thm:main1}(2) to $\piab(\CX)$, we obtain a type
  $q(x)\in S_{\piab(\CX)}(A)$ that is $\Gamab$-dense in $\Gab$. Let
  $p\in S_\CX(A)$ be a type on $\CX$ with $\piab(p)=q$.  By
  Theorem~\ref{thm:type-dense}, the type $p$ is $\Gamma$-dense in $G$,
  hence $\CX$ is $\Gamma$-dense in $G$ as well.
\end{proof}

\section{Interpreting the results as Hausdorff limits}

\subsection{Hausdorff limits}

We first recall some definitions.

Let $(M,d)$ be a compact metric space, and $X_1,X_2\subseteq X$.  The
Hausdorff distance $d_H(X_1,X_2)$ between $X_1$ and $X_2$ is defined
as follows: First, for $x\in M$, we let
$d(x,X_i)=\inf_{y\in X_i}d(x,y)$. Next,
\[d_H(X_1,X_2)=\max\{\sup_{x\in X_1} d(x,X_2),\, \sup_{x\in X_2}d(x,X_1)\}.\]

An equivalent definition is given by:
\[d_H(X_1,X_2)=\inf\{r\geq 0: \forall x_i\in X_i, i=1,2,\,\,\,
d(x_1,X_2),d(x_2,X_1)\leq r \}.\]

\begin{rem}\label{rem:hlim} For $X_1,X_2 \subseteq M$ we have
  $d_H(X_1,X_2)=0$ if and only if $\cl(X_1)=\cl(X_2)$.
\end{rem}
Denoting by $\CK(M)$ the set of all compact subsets of $M$, it is
known that the restriction of $d_H$ to $\CK(M)$ makes it into a
compact metric space.

The topology induced by $d_H$ on $\CK(M)$ does not depend on the
metric $d$ but only on the topology of $M$. It coincides with the
Vietoris topology.

Given a family $\CF \subseteq \CK(M)$, a set $Y\in \CK(M)$ is \emph{a
  Hausdorff limit} of $\CF$ if for every $\varepsilon>0$ there is
$F\in \CF$ with $d_H(Y,F)<\varepsilon$. Using Remark~\ref{rem:hlim},
we extend this definition to a family $\CF\subseteq \CP(M)$ of
arbitrary subsets of $M$ by saying that $Y\in \CK(M)$ is a Hausdorff
limit of $\CF$ if it is a Hausdorff limit of the family
$\{ \cl(F) \colon F\in \CF \}$.

\subsubsection{Limits at infinity}\label{sec:limits-at-infinity}

We denote by $I_\infty$ the interval $(0,+\infty) \subseteq \RR$. We
abreviate ``for all sufficiently large $t$'' by ``$t\gg 0$''.

We define:
\begin{defn}\label{def:converges at infty}
  Let $\CF =\{ F_t \colon t\in I_\infty\}$ be a family of subsets of
  $M$.
  \begin{enumerate}
  \item A set $Y\in \CK(M)$ is \emph{a Hausdorff limit at $\infty$ of
      the family $\CF$} if for all $\epsilon >0$ and $r>0$ there is
    $t> r$ with $d_H(Y, F_t) < \varepsilon$.
  \item We say that \emph{the family $\CF$ converges to a set
      $Y\in \CK(M)$ at $\infty$} if $Y$ is the unique Hausdorff limit
    of $\CF$ at $\infty$. In this case, since $\CK(M)$ is compact, $Y$
    is \textbf{the} limit of $\CF$ as $t$ goes to $\infty$, namely for
    any $\epsilon> 0$ there is $R\in \RR$ such that for all $t>R$,
    $d_H(Y,F_t)< \varepsilon$.
  \end{enumerate}

\end{defn}

\subsection{Haudorff Limits via the standard part map.}\label{sec:haudorff-limits-via}

We fix a compact set $M\subseteq \RR^n$ with the metric $d$ induced by
the Euclidean metric of $\RR^n$, and we view both $M$ and $d$ as
definable in $\Rfull$.

Since $M$ is compact, $M^\sharp\subseteq \CO^n$, and we denote by
$\st_M$ the restriction of the standrad part map
$\st\colon \CO^n\to \RR^n$ to $M^\sharp$.  It is not hard to see that
$\st_{M}\colon M^\sharp \to M$ maps $\alpha\in M^\sharp$ to the unique
$a\in M$ such that $\alpha\in U^\sharp$ for every neighborhood $U$ of
$a$.

\medskip

Let $\CF=\{ F_t \colon t\in T\}$ be a family of subsets of $M$ indexed
by a set $T\subseteq \RR^m$.  We can view this family also as the
family of fibers of the set
$F=\{ (x,t)\in X\times T \colon x\in F_t\}$, with respect to the
second projection, and hence as a family definable in $\Rfull$.  Thus
for $\tau\in T^\sharp$ we also have a ``non-standard'' fiber
$F_\tau^\sharp=\{ x\in M^\sharp : (x,\tau)\in F^\sharp \}$.  Using
\cite{narens}*{Theorem 4.4} we obtain:

\begin{fact}\label{fact:narens}
  In the above setting, a set $Y\in \CK(M)$ is a Hausdorff limit of
  the family $\CF=\{ F_t \colon t\in T\}$ if and only if there is
  $\tau\in T^\sharp$ such that $Y=\st_M(F^\sharp_{\tau})$.
\end{fact}

Using the above, we conclude:

\begin{lem}\label{claim:hinfty}
  For a family $\CF =\{ F_t \colon t\in I_\infty\}$ of subsets of $M$
  indexed by $I_\infty$, a set $Y\in \CK(M)$ is a Hausdorff limit at
  $\infty$ of $\CF$ if and only if there is $\tau \in \mfR$ with
  $\tau > \RR$, such that $Y=\st_M(F^\sharp_\tau)$.
\end{lem}
\begin{proof}
  For $r\in \RR^{>0}$ let $I_{>r}$ be the interval
  $(r,+\infty)\subseteq \RR$ and $\CF_r$ be the family
  $\CF_r=\{ F_t \colon t\in I_{>r}\}$.

  It is easy to see that a set $Y\in \CK(M)$ is a Hausdorf limit at
  $\infty$ of the family $\CF$ if and only if for every
  $r\in \RR^{>0}$ the set $Y$ is a Hausdorf limit of the family
  $\CF_r$.

  By Fact~\ref{fact:narens}, the latter condition is equivalent to the
  following: for every $r\in \RR^{>0}$, there is
  $\tau_r \in I^\sharp_{>r}$ with $Y=\st_M(F^\sharp_{\tau_r})$.

  Thus the conclusion of the lemma can be restated as follows:

  \medskip
  \begin{minipage}[t]{0.9\linewidth}
    For a set $Y\in \CK(M)$ the following are equivalent:
    \begin{enumerate}[(a), leftmargin = 2pc]
    \item For every $r\in \RR^{>0}$, there is
      $\tau_r \in I^\sharp_{>r}$ with $Y=\st_M(F^\sharp_{\tau_r})$.
    \item There is $\tau \in\mfR$ with $\tau > \RR$ such that
      $Y=\st_M(F^\sharp_\tau)$.
    \end{enumerate}
  \end{minipage}\\

  The direction $(b)\Rightarrow (a)$ is obvious, and the opposite
  direction follows from the $|\RR|^+$-saturation of $\fRfull$.

\end{proof}

\subsection{Hausdorff limits  in $\mathbf{G/H}$}
Let $G$ be a connected Lie group and $H\subseteq G$ a closed subgroup
such that the space of the left cosets $N=G/H$ is compact, with
respect to the quotient topology.  We denote by $\pi\colon G\to N$ the
quotient map.  Using Whitney embedding theorem we embed $G$ into some
$\RR^m$ and $N$ into some $\RR^n$ as closed subsets, and view $G$, $N$
and $\pi$ as definable in $\Rfull$.

\medskip

Given a family $\CF=\{ F_t \colon t\in T\}$ of subsets of $G$, we let
\[\pi(\CF)=\{\pi(F_t) \colon t\in T \}\]
be the corresponding family of subsets of $N$.

\begin{prop}\label{prop:hlimpi}
  Let $\CF=\{ F_t\colon t\in I_\infty \}$ be a family of subsets of
  $G$.  A set $Y\in \CK(N)$ is a Hausdorf limit at $\infty$ of the
  family $\pi(\CF)$ if and only if there is $\tau\in \mfR$ with
  $\tau > \RR$ such that
  $Y=\pi\bigl( \st_G( F^\sharp_\tau \ccdot H^\sharp )\bigr)$.
\end{prop}
\begin{proof}
  For $\tau \in \fRfull$, by Claim~\ref{claim:hinfty}, it is
  sufficient to show that
  \[ \st_{N} (\pi^\sharp(F_\tau^\sharp))=\pi\bigl( \st_G(
    F^\sharp_\tau \ccdot H^\sharp )\bigr). \] For
  $\alpha\in G^\sharp$, we will show that
  $ \st_{N} (\pi^\sharp(\alpha))=\pi\bigl( \st_G( \alpha \ccdot
  H^\sharp )\bigr)$.  That is clearly enough.

  Since $G/H$ is compact, there is $\beta \in \CO^m\cap G^\sharp$ with
  $\beta\in \alpha\ccdot H^\sharp$.  Let $b=\st_G(\beta)$.  Since $H$
  is a closed subgroup, the set $b\ccdot H$ is closed and, by
  Lemma~\ref{lem:st-closed}(2), we have
  \[ b\ccdot H = \st_G(b\ccdot H^\sharp)=\st_G(\beta\ccdot
    H^\sharp). \] Since $\pi^\sharp$ is invariant under the action of
  $H^\sharp$ on the right we also have
  $\pi^\sharp(\alpha) =\pi^\sharp(\beta)$ and we are left to show
  \[ \st_{N} (\pi^\sharp(\beta))=\pi(b). \] Since $\pi$ is continuous,
  the latter follows from Fact~\ref{fact:st-com}.
\end{proof}

\subsection{Hausdorff limits in nilmanifolds}
We go back to our o-minimal structure $\Rom$ and fix a unipotent group
$G$.

For a lattice $\Gamma\subseteq G$, we use $\pi_\Gamma$ to denote the
projection $\pi_\Gamma\colon G\to G/\Gamma$. When no confusion arises,
we omit the subscript $\Gamma$. Also, whenever $\Gamma_0\subseteq \Gamma$
is a subgroup of finite index, we let $\pi_0\colon G\to G/\Gamma_0$
denote the natural projection.

Given an $\Rom$-definable family $\CF=\{F_t:t\in I_\infty\}$, for a
lattice $\Gamma\subseteq G$ we consider the possible Hausdorff limits at
$\infty$ of the family $\pi(\CF) \subseteq G/\Gamma$. Notice that if $\CF$
is a constant family $F_t=F$ then the only Hausdorff limit at $\infty$
is the closure of $\pi(F)$ and this case was handled in
\cite{nilpotent}.

\begin{sample}\label{sample:hausd}

  \begin{enumerate}
  \item Consider first $G=(\RR^2,+)$ and $\Gamma=\ZZ^2$.

    Let $L_0$ be the line
    $L_0= \{ (x,0) \in \RR^2 \colon x \in \RR\}$, and $\CF_1$ be the
    family of $L_0$-translates:
    $\CF_1=\{ L_0+(0,t) \colon t\in I_\infty \}$.  It is not hard to
    see that the Hausdorff limits at $\infty$ of $\pi(\CF_1)$ are
    exactly the sets $\pi(L_0+g)$ for $g\in G$.

    Let $\CF_2=\{ L_t \colon t\in I_\infty \}$ be the family of lines
    in $G$ where $L_t $ is the line $L_t=\{ (x,y)\in G \colon y=tx\}$.
    It is not hard to see that the only Hausdorff limit at $\infty$ is
    the whole $G/\Gamma$.

  \item Assume now that $G=\RR$, $\Gamma=\ZZ$, and let
    $\CF=\{t+[0,2]:t\in I_\infty\}$.  The family $\pi(\CF)$ is the
    constant family $\pi(F_t)=\RR/\ZZ$, and hence this is the only
    Hausdorff limit at $\infty$.  However, for any lattice
    $\Gamma_0 \subseteq \ZZ$ with $|\ZZ:\Gamma_0|\geq 3$, the
    Hausdorff limits at $\infty$ of $\pi_0(\CF)$ are the sets of the
    form $g+\pi_0([0,2])$, for $g\in \RR/\Gamma_0$, and none of these
    equals $\RR/\Gamma_0$.
  \end{enumerate}
\end{sample}

\begin{defn}
  Let $\CF=\{F_t:t\in I_\infty\}$ be an $\Rom$-definable family of
  subsets of a unipotent group $G$, and $\Gamma\subseteq G$ be a lattice.

  We say that the family $\pi(\CF)$ \emph{converges strongly} to
  $G/\Gamma$ at $\infty$ if $\pi_0(\CF)$ converges to $G/\Gamma_0$ at
  $\infty$, for any subgroup $\Gamma_0\subseteq \Gamma$ of finite
  index, as in Definition \ref{def:converges at infty}.

\end{defn}
The next observation immediately follows from
Proposition~\ref{prop:hlimpi}:
\begin{prop}\label{prop: dense and converge} Let
  $\CF=\{F_t:t\in I_\infty\}$ be an $\Rom$-definable family of subsets
  of a unipotent group $G$, and $\Gamma\subseteq G$ be a lattice. Then:
  \begin{enumerate}
  \item $G/\Gamma$ is a Hausdorff limit of $\pi(\CF)$
    $\Leftrightarrow$ there exists $\tau\in \mfR$, $\tau >\mathbb R$,
    such that $F_\tau^\sharp$ is $\Gamma$-dense in $G$.
  \item $\pi(\CF)$ converges to $G/\Gamma$ $\Leftrightarrow$ for all
    $\tau\in \mfR$ with $\tau >\mathbb R$, $F_\tau^\sharp$ is
    $\Gamma$-dense in $G$.
  \item $\pi(\CF)$ converges strongly to $G/\Gamma$ $\Leftrightarrow$
    for all $\tau\in \mfR$, with $\tau >\mathbb R$, $F_\tau^\sharp$ is
    strongly $\Gamma$-dense in $G$.
  \end{enumerate}
\end{prop}

For the next result, recall the notation at the beginning of Section
\ref{sec:abelinization}, regarding the abelianization of $G$.

\begin{cor}\label{cor:abel}
  Let $G$ be a unipotent group and let $\CF=\{F_t:t\in I_\infty\}$ be
  an $\Lom$-definable family of subsets of $G$. We denote by
  $\CF_\mathrm{ab}$ the family $\piab(\CF)$ of subsets of $\Gab$. For
  a lattice $\Gamma\subseteq G$, we let $\pi:G\to G/\Gamma$, and
  $\pi^*:\Gab\to \Gab/\Gamab$ be the quotient maps.

  Then, the family $\pi(\CF)$ converges strongly to $G/\Gamma$ at
  $\infty$ if and only if $\pi^*(\CF)$ converges strongly to
  $\Gab/\Gamab$ at $\infty$.
\end{cor}
\begin{proof} By Proposition \ref{prop: dense and converge},
  $\pi(\CF)$ converges strongly to $G/\Gamma$ at $\infty$ if and only
  if for all $\tau>\RR$ in $\mfR$, $F^\sharp_{\tau}$ is strongly
  $\Gamma$-dense in $G$. By Corollary \ref{cor:xiffxab}, this is
  equivalent to $\piab(F^\sharp_\tau)$ being strongly $\Gamab$-dense
  in $\Gab$, for all $\tau>\RR$, which again, by Proposition
  \ref{prop: dense and converge}, is equivalent to $\piab(\CF)$
  strongly converging to $\Gab/\Gamab$ at $\infty$.

\end{proof}

Before the next theorem we observe:
\begin{lem}\label{lem: eventual Lmax} Let $G$ be a unipotent group
  and let $\CF=\{F_t:t\in I_\infty\}$ be an $\Rom$-definable family of
  subsets of $G$. Then, for all $\tau,\tau'\in \mfR$, with
  $\tau,\tau'>\mathbb R$,
  $\SLm(F_{\tau}^\sharp)=\SLm(F_{\tau'}^\sharp)$.
\end{lem}
\begin{proof} We use the fact that $\tau$ and $\tau'$ have the same
  $\Lom$-type over $\RR$.

  Clearly, it is enough to show
  $\SL(F^\sharp_\tau)=\SL(F^\sharp_{\tau'})$, and, by symmetry, it is
  sufficient to show
  $\SL(F^\sharp_\tau)\subseteq \SL(F^\sharp_{\tau'})$.

  Let $L\in \SL(F^\sharp_\tau)$. We choose $\alpha\in \dcl(\tau)$ such
  that the coset $\alpha L$ is a nearest co-commutative coset to some
  type on $F_{\tau}^\sharp$.  Let $a(t)$ be an $\Rom$-definable
  function with $a(\tau)=\alpha$,

  By saturation of $\fRom$, the coset $a(\tau')L$ is a nearest
  co-commutative coset to some type on $F_{\tau'}^\sharp$.
\end{proof}

The above lemma justifies the following definition
\begin{defn}
  For an $\Rom$-definable family $\CF=\{F_t:t\in I_\infty\}$ of a
  unipotent group $G$, we denote by $\SLm(\CF)$ the finite set of
  co-commutative subgroups $\SLm(F_{\tau}^\sharp)$, for some (any)
  $\tau>\RR$.
\end{defn}

The next theorem is one of our main results.
\begin{thm}\label{thm:main-hlin1}
  Let $G$ be a unipotent group, $\CF=\{F_t:t\in I_\infty\}$ an
  $\Rom$-definable family of subsets of $G$.

  For every lattice $\Gamma\subseteq G$ we have:
  \begin{enumerate}[resume]
  \item $L^\Gamma=G$ for some $L\in \SLm(\CF)$ if and only if
    $\pi(\CF)$ converges strongly to $G/\Gamma$ at $\infty$.

  \item $L^\Gamma\neq G$ for all $L\in \SLm(\CF)$ if and only if there
    exists a subgroup $\Gamma_0\subseteq \Gamma$ of finite index such
    that all Hausdorff limits at $\infty$ of $\pi_0(\CF)$ are proper
    subsets of $G/\Gamma_0$.
  \end{enumerate}
\end{thm}

Note that the condition given in (2) is formally stronger than the
negation of the condition in (1), thus both need to be proved
separately.

\begin{proof} (1) By Proposition~\ref{prop: dense and converge},
  $\pi(\CF)$ converges strongly to $G/\Gamma$ if and only if for all
  $\tau>\mathbb R$, $F_{\tau}^\sharp$ is strongly $\Gamma$-dense in
  $G$, which by Theorem~\ref{thm:main1}, is equivalent to $L^\Gamma=G$
  for some $L\in \SLm(\CF)$.

  (2) Assume that for all $L\in \SLm(\CF)$ we have $L^ \Gamma\neq G$,
  and for contradiction assume that for every subgroup
  $\Gamma_0\subseteq \Gamma$ of finite index, $G/\Gamma_0$ is one of
  the Hausdorff limits at $\infty$, of the family $\pi_0(\CF)$.
  Equivalently, it follows from Proposition~\ref{prop:hlimpi} and
  Lemma \ref{lem:Gamma-dense}, that for every subgroup
  $\Gamma_0\subseteq \Gamma$ of finite index, there is a
  $\tau'\in \mfR$ with $\tau'>\RR$, such that
  $G=\st(F^\sharp_{\tau'}\ccdot \Gamma_0^\sharp)$.

  \begin{claim} There exists $\tau^*\in \mfR$ with $\tau^*>\RR$ such
    that for every $\Gamma_0\subseteq \Gamma$ of finite index,
    $\st(F^{\sharp}_{\tau*}\ccdot\Gamma_0^\sharp)=G$.

  \end{claim}
  \begin{proof}[Proof of the claim.]
    For every $g\in G=G(\RR)$, $r\in \RR^{>0}$ and
    $\Gamma_0\subseteq \Gamma$ of finite index, we consider the following
    formula
    $\phi_{g,r,\Gamma_0}(t)$:
    \[t>r\,\&\, g\in B_{1/r}(e)\ccdot
    F_t\ccdot \Gamma_0.\]

    We let $p(t)$ be the type consisting of all $\phi_{g,r,\Gamma_0}$,
    as $g,r,\Gamma_0$ vary over all $g\in G$, $r\in \RR^{>0}$ and
    $\Gamma_0\subseteq \Gamma$ of finite index, respectively. We claim that
    $p$ is finitely consistent. Indeed, given finitely many subgroups
    of $\Gamma$ of finite index, let $\Gamma_1$ be their
    intersection. Clearly $\Gamma_1$ has finite index in $\Gamma$.  By
    our assumption, there is $\tau>\RR$ such that
    $G(\RR)\subseteq \mu\ccdot F^\sharp_{\tau}\ccdot\Gamma_1^\sharp$, which
    implies that for any $g_1,\ldots, g_k\in G$ and
    $r_1,\ldots, r_k\in \RR$, $\phi_{g_i,r_i, \Gamma_1}(\tau)$
    holds. It follows that $p(t)$ is finitely consistent, so by the
    saturation of $\fRfull$, there exists $\tau^*\in \mfR$ realizing
    $p(t)$.

    Now, given $\Gamma_0\subseteq \Gamma$ of finite index, and $g\in G$, we
    have
    $g\in B_\epsilon(e)^\sharp \ccdot F_{\tau*}^\sharp
    \ccdot\Gamma_0^\sharp$, for all $\epsilon\in \RR^{>0}$. Using
    saturation again, it follows that
    $g\in \mu\ccdot F_{\tau*}^\sharp\ccdot \Gamma_0^\sharp$, and hence
    $G=\st(F_{\tau*}^\sharp \ccdot \Gamma_0^\sharp)$, proving the
    claim.
  \end{proof}

  For $\tau^*$ as in the above claim, $F^\sharp_{\tau*}$ is strongly
  $\Gamma$-dense in $G$ and therefore, by Theorem~\ref{thm:main1},
  there is $L\in \SLm(\CF)$ with $L^\Gamma=G$, contradiction.

  The opposite implication of (2) follows from (1).

\end{proof}

\subsection{The abelian case}
In the unipotent case, Theorem~\ref{thm:main-hlin1} tells us when the
family $\pi(\CF)$ converges (strongly) to $G/\Gamma$. In the abelian
case we can say more about the possible Hausdorff limits of
$\pi(\CF)$, due to the following theorem.

\begin{thm}\label{thm:main-hlin-ab}
  Let $G=(\RR^m,+)$ and let $\CF=\{F_t:t\in I_\infty\}$ be an
  $\Rom$-definable family of subsets of $G$.

  For every $r\in \RR^{>0}$, there are subspaces
  $L_1,\ldots, L_k\subseteq G$ (possibly with repetitions and with $k$
  depending on $r$), with $\SLm(\CF)\subseteq \{L_1,\ldots, L_k\}$, and
  there are $\Rom$-definable functions
  $a_1(t),\ldots, a_{k}(t)\colon I_\infty \to G$, such that
  \begin{enumerate}

  \item For some $\tau\in \mfR, \tau>\RR$, each $a_i(\tau)+L_i$ is a
    nearest coset to some type $p\in S_{F_\tau^\sharp}(\tau)$.
  \item For $t \gg 0$,
    \[ F_t\subseteq \bB_r +\bigcup_{i=1}^{k} a_i(t) + L_i.\]
  \item Let $\Gamma\subseteq G$ be a lattice and $\pi:G\to G/\Gamma$ the
    projection. For every $s\in \RR^{>0}$, there exists $t_s>0$ such
    that for all $t>t_s$, and all $i=1,\ldots, k$.
\[\pi(a_i(t)+ L_i)\subseteq \pi(\bB_s) +\pi(F_t).\]
\end{enumerate}
\end{thm}

\begin{proof} Let $\tau\in I_\infty^\sharp$ with $\tau > \RR$, and let
  $A=\dcl(\tau)$.

  Fix $r\in \RR^{>0}$. Let $L_1,\ldots,L_{k}\in \SL_A(F_\tau^\sharp)$
  be subgroups and $\alpha_1,\dotsc,\alpha_{k}\in \dcl(A)$ be as in
  Theorem~\ref{thm:main-section3}.  Thus we have
\[ F_\tau^\sharp\subseteq  \bB_r^\sharp + \bigcup_{i=1}^{k} \alpha_i +
L_i^\sharp,\,\,\,\, \mbox{ with } \SLm(F_\tau^\sharp)\subseteq
\{L_1,\ldots,L_k\},\] and, by Theorem~\ref{thm:main types}, for every
$i=1,\dotsc,k$, and $s>0$ we also have
\[\alpha_i+L_i^\sharp \subseteq \bB_s^\sharp +F_\tau^\sharp
+\Gamma^\sharp.\] Since each $\alpha_i\in \dcl(\tau)$, for
$i=1,\dotsc,k$, we choose $\Rom$-definable functions
$a_i(t)\colon \RR\to G$, such that $\alpha_i=a_i(\tau)$.  We have
\[
  F_\tau^\sharp\subseteq \bigcup_{i=1}^{k} \bB_r^\sharp + a_i(\tau) +
  L_i^\sharp.
\]

Since the $\Lom$-type of $\tau$ over $\RR$ is implied by
$\{ x> r \colon r\in \RR\}$ and the above inclusion can be expressed
by an $\Lom$-formula over $\tau$, we obtain that for $t\gg 0$ ,
\[F_t\subseteq \bigcup_{i=1}^{k} \bB_r + a_i(t) + L_i.\] This proves (1)
and (2).

\medskip Assume (3) fails. Then, there is a lattice
$\Gamma \subseteq G$ and $s\in \RR^{>0}$, such that for some
$i_0=1,\dotsc,k$, the set
\[\{ t\in I_\infty \colon \pi(a_i(t)+ L_i) \nsubseteq \pi(\bB_s)
+\pi(F_t)\}\] is unbounded in $\RR$.  Without loss of generality, we
assume $i_0=1$.

Using saturation of $\fRfull$, we can find $\tau'\in \mfR$ with
$\tau'>\RR$, such that
$\pi (a_1(\tau')+ L_1^\sharp) \nsubseteq \pi(\bB_s^\sharp)
+\pi(F_{\tau'}^\sharp)$, so in particular,
\[ a_1(\tau')+ L_1^\sharp \nsubseteq \bB_s^\sharp +F_{\tau'}^\sharp
  +\Gamma^\sharp.\] Since $\tau$ and $\tau'$ realize the same
$\Lom$-type over $\RR$, and $a_1(\tau)+L_1^\sharp$ is the nearest
coset to some type on $F^\sharp_\tau$, the coset
$a_1(\tau')+L_1^\sharp$ is the nearest coset to a type, call it $p$,
on $F^\sharp_{\tau'}$. However, using Theorem~\ref{thm:main types}, we
have
\[ a_1(\tau')+L_1^\sharp\subseteq \mu+a_1(\tau')+L_1^\sharp=\mu+p(\mfR)+\Gamma^\sharp\subseteq \bB_s^\sharp +F_{\tau'}^\sharp +\Gamma^\sharp,\] contradiction.

\end{proof}

We can now show that, in the abelian case, every Hausdorff limit of
$\pi(\CF)$ at $\infty$ is trapped between a finite union of cosets and
a ``thickening'' of it.
\begin{cor}\label{cor:abelian-hl}
  Let $G=(\RR^m,+)$ and let $\CF=\{F_t:t\in I_\infty\}$ be an
  $\Rom$-definable family of subsets of $G$.

  For every $r\in \RR^{>0}$, there are subspaces
  $L_1,\ldots, L_k\subseteq G$ (possibly with repetitions and with $k$
  depending on $r$), with $\SLm(\CF)\subseteq \{L_1,\ldots, L_k\}$, such
  that for any lattice $\Gamma\subseteq G$, and any Hausdorff limit
  $X$ of the family $\pi(\CF)$ at $\infty$, there are
  $g_1,\dotsc,g_k\in G$ with
  \[ \pi\Bigl(\bigcup_{i=1}^k( g_i +L_i)\Bigr) \subseteq X
    \subseteq\pi(\bB_r)+\pi\Bigl(\bigcup_{i=1}^k( g_i +L_i)\Bigr). \]
\end{cor}

\begin{proof}
  By Proposition \ref{prop:hlimpi},
  $X=\pi(\st(F_\tau^\sharp +\Gamma^\sharp))$, for some $\tau>\RR$. We
  now apply Theorem \ref{thm:main-hlin-ab} and obtain linear subspaces
  $L_1,\ldots, L_k\subseteq \RR^m$ and definable functions
  $a_1(t),\ldots, a_k(t)\colon I_\infty\to G$ as in the theorem. Given
  a lattice $\Gamma\subseteq G$, clause (2) of that theorem implies that
\[F_\tau^\sharp+\Gamma^\sharp\subseteq \bB_r^\sharp+\bigcup_{i=1}^k a_i(\tau)+L_i^\sharp+\Gamma^\sharp,\] while  (3) (applied for all $r\in \RR^{>0}$) implies
\[\bigcup_{i=1}^k a_i(\tau)+L_i^\sharp+\Gamma^\sharp \subseteq\mu+ F_\tau^\sharp +\Gamma^\sharp.\]

Using Lemma~\ref{lem:fst-cyl}, for
$g_i\in \st(a_i(\tau)+\Gamma^\sharp)$, we obtain
\[\bigcup_{i=1}^k g_i+L_i^\Gamma+\Gamma
=\bigcup_{i=1}^k \st(a_i(\tau)+L_i^\sharp+\Gamma^\sharp)
\subseteq\st(F_\tau^\sharp +\Gamma^\sharp)\subseteq \bB_r +\bigcup_{i=1}^k
g_i+L_i^\Gamma+\Gamma.\]

Applying $\pi$ the result follows.\end{proof}

\section{Polynomial dilations}\label{sec: polynomial}

Let $G\subseteq \GL(n,\mathbb R)$ be a unipotent group. In several articles
(\cite{KSS}, \cite{fish}) a particular type of families of subsets of
$G$, given by dilations of an initial curve, was considered in the
unipotent setting. We first make some definitions.
\subsection{The setting}\label{sec: setting}

\begin{defn} \emph{A polynomial $m\times k$ matrix $M_t$ (over
    $\mathbb R)$} is a matrix $M_t\in M_{m\times k}(\RR[t])$, namely,
  a matrix all of whose entries are polynomials in $\mathbb R[t]$.  It
  can be written as $\sum_{j=0}^d t^j A_j$, with each $A_j$ is an
  $m\times k$ matrix over $\RR$.
\end{defn}

Let $G$ be a unipotent group of dimension $m$. We identify the
underlying vector space of its Lie algebra $\mfg$ with $\RR^m$, and we
let $\exp:\mfg\to G$ be the exponential map (which is a polynomial
bijection with a polynomial inverse, see \cite{nilpotent-book}).  For
a polynomial $m\times k$ matrix $M_t$ we consider the family of
``dilations'' $\rho_t\colon \RR^k\to G$ given by $x\mapsto \exp(M_tx)$
and for a set $X\subseteq \mathbb R^k$, the family
$\{\rho_t(X):t\in I_\infty\}$ of subsets of $G$.

In \cite{KSS}, the authors start with a measure $\nu$ on $\mfg$ given
as the pushforward of the Lebesgue measure on $(0,1)$ via a real
analytic map $\phi:(0,1)\to \mfg$, then ``dilate'' the measure $\nu$
using multiplication by a polynomial $m\times m$ matrix $M_t$, and
consider the limit of the measures as $t\to \infty$. In \cite[Theorem
1.1]{KSS} the authors prove that under some assumptions on
$\Phi=\mathrm{Image}(\phi)$, given in terms of kernels of particular
characters of $G$, the associated family of measures $\mu_t$, on
$G/\Gamma$ is ``equidistributed'', roughly saying that for any Borel
$D\subseteq G/\Gamma$, the family $\mu_t(D)$ converges to the canonical
Haar measure of $D$. Translated to the topological language this would
imply (but a-priori might not be equivalent) that the\ family
$\{\pi_\Gamma(\rho_t(\Phi))\colon t\in I_\infty\}$ converges to
$G/\Gamma$ at $\infty$.

In the current section we extend the topological corollary by
replacing the one-dimensional set $\Phi$ with an $\Rom$-definable set
$X\subseteq \mathfrak \RR^k$ of arbitrary dimension.  Using
Theorem~\ref{thm:main-hlin1}, one can formulate conditions, similar to
those in \cite{KSS}, as to when the family
$\{\pi(\exp(M_t\ccdot X)) \colon t\in I_\infty\}$ converges strongly
to $G/\Gamma$ at $\infty$.  Instead, we consider the special case,
when the polynomial matrix $M_t$ does not have a constant term, and
describe first in full, in the abelian case, all possible Hausdroff
limits at $\infty$ of the family.

\subsection{Polynomial dilations in vector spaces}\label{sec:dial-vect-spac}
We fix some notations.

\begin{ntn}
  \begin{enumerate}
  \item We call a map $a(t):\RR\to \RR^m$ \emph{a polynomial map} if
    $a(t)=\sum_{i=0}^d t^ia_i$ for some $a_0,\dotsc,a_d\in \RR^m$. We
    say that the polynomial map is \emph{proper} if $a_0$=0.
  \item By \emph{a polynomial family of dilations} we mean a family of
    maps $\{\rho_t\colon \RR^k\to \RR^m \colon t\in I\}$ with
    $I\subseteq \RR$, such that for some polynomial $m\times k$-matrix
    $M_t$, for all $t\in I$ and $x\in \RR^k$, we have
    $\rho_t(x)=M_tx$.  We say that the family is \emph{proper} if the
    constant term of $M_t$ is the zero $m\times k$-matrix.
  \item By a \emph{a polynomial family of cosets} we mean a family
    $\{ a(t)+L \colon t\in I\}$, where $I\subseteq \RR$,
    $L\subseteq \RR^m$ is a subspace, and $a(t)\colon \RR \to \RR^m$ a
    polynomial map, We say that the family is \emph{proper} if $a(t)$
    is a proper polynomial map.
  \item By \emph{a polynomial family of multi-cosets} we mean a family
    of the form $\{ \bigcup_{j=1}^n(a_j(t)+L_j) \colon t\in I\}$,
    where $I\subseteq \RR$, each $a_j(t)\colon \RR\to \RR^m$ is a
    polynomial map and each $L_j\subseteq \RR^m$ is a subspace.  We
    say that the polynomial family of multi-cosets is \emph{proper} if
    each $a_j(t)$ is proper.
  \end{enumerate}
\end{ntn}
We refer the readers to \cite[Remark 2.5]{KSS} for an example
explaining the reason to work with proper polynomial dilations instead
of general ones.

  \begin{rem}\label{rem:multimax}
    Let $\CM=\{ \bigcup_{j=1}^n(a_j(t)+L_j) \colon t\in I\}$ be a
    polynomial family of multi-cosets. It is not hard to see, using
    Lemma~\ref{lem:fst-cyl}(2), that $\SLm(\CM)$ is exactly the set of
    maximal (by inclusion) subspaces in $\{L_1,\dotsc,L_n\}$.
  \end{rem}

  Our main result is:

\begin{thm}\label{thm:main-dilation}
  Let $\{ \rho_t\colon \RR^k\to \RR^m\colon t\in I_\infty\}$ be a
  proper family of polynomial dilations, $X\subseteq \RR^k$ an
  $\Rom$-definable set, and $\CF=\{ \rho_t(X) \colon t\in I_\infty\}$.

  Then, there is a proper polynomial family of
  multi-cosets
  \[\CM=\Bigl\{ \bigcup_{j=1}^n(p_j(t)+L_j) \colon t\in
  I_\infty\Bigr\}\] such that
  \begin{enumerate}
  \item $\rho_t(X)\subseteq \bigcup_{j=1}^n(p_j(t)+L_j) $ for all
    $t\in \RR$;
  \item $\SLm(\CF)=\SLm (\CM)$;
  \item for any lattice $\Gamma\subseteq \RR^m$, the families
    $\pi(\CF)$ and $\pi(\CM)$ have the same Hausdorff limits at
    $\infty$.
  \end{enumerate}
\end{thm}

\begin{proof}
  Let $A_1,\dotsc, A_d$ be $m\times k$-matrices such that
  $\rho_t(x)=\sum_{i=1}^d t^i A_i x$.

  We fix $\tau\in \mfR$ with $\tau>\RR$, and let
  $\CY= \rho_\tau(X^\sharp)$.

  \smallskip Applying Theorem~\ref{thm:main-hlin-ab} with $r=1$, we
  obtain subspaces $L_1,\dotsc L_n\subseteq \RR^m$ and
  $\Rom$-definable maps
  $a_1(t), \dotsc, a_n(t)\colon I_\infty \to \RR^m$ such that
  \begin{equation}
    \label{eq:5}
    \begin{gathered}
      \text{each } a_j(\tau)+L_j \text{ is a nearest coset to some
        type in }S_\CY(\tau),
      \text{ and } \\
      \rho_t(X) \subseteq \bB_1+\bigcup_{j=1}^n \bigl(a_j(t)+L_j
      \bigr) \text{ for } t\gg 0.
    \end{gathered}
  \end{equation}

  We pick $a_i(t)$ and $ L_i$, $ j=1,\dotsc,n$, as above, with the
  minimal possible $n$.

  We claim that $\SLm(\CF)$ is exactly the set of maximal elements of
  $\{L_1,\dotsc,L_n\}$. Indeed, assume $L\in \SLm(\CF)$ and
  $p\in S_{\CY}(\tau)$ is such that $L_p=L$. Then, for some
  $j=1,\ldots, n$, $p\vdash a_j(\tau)+L_j+\bB_1$. It follows from
  Lemma~\ref{lem:properties} that $L\subseteq L_j$, so by maximality
  $L=L_j$. Since every $L_i$ is contained in some $L\in \SLm(\CF)$,
  the claim follows.

  \medskip

  For each $j=1,\dotsc,n$, let $L_j^\perp\subseteq\RR^m$ be the
  orthogonal complement to $L_j$ with respect to the standard inner
  product on $\RR^m$, and $\pi_j^\perp\colon \RR^m\to L_j^\perp$ be
  the corresponding projection, whose kernel is $L_j$.

  Replacing each $a_j(t)$ with $\pi_j^\perp(a_j(t))$, if needed, we
  assume $a_j(t)\in L_j^\perp$ for all $t\in I_\infty$.

  For $j=1,\dotsc,n$, let
  \[ X_j =\{ x\in X \colon \text{ for } t\gg 0, \,\, \rho_t (x) \in
    \bB_1+ a_j(t)+ L_j \}. \]

  Notice that each $X_j$ is $\Rom$-definable,
  \[ X=\bigcup_{j=1}^n X_j, \] and, by the minimality of $n$, each
  $X_j$ is non-empty.

  \smallskip
  \noindent\textbf{Claim A.} {\itshape For each $j=1,\dotsc,n$, and
    every $i=1,\dotsc,d$, the set $A_iX_j$ is contained in a single
    coset of $L_j$.}
  \begin{proof}[Proof of Claim A]
    We fix $j\in\{1,\dotsc,n\}$.

    Clearly it is sufficient to show that for any $i=1,\dotsc,d$,
    every $\RR$-linear function $F\colon \RR^m\to \RR$ that vanishes
    on $L_j$ is constant on $A_iX_j$.

    Assume not, and for some $\RR$-linear function
    $F\colon \RR^m\to \RR$ vanishing on $L_j$ there are
    $r \in \{1,\dotsc,d\}$ and $b_1, b_2\in X_j$ with
    $F(A_rb_1)\neq F(A_rb_2)$.

    Consider the map $q(t)= F(\rho_t (b_1)-\rho_t(b_2))$.

    Since $\rho_t (b_1)-\rho_t(b_2)= \sum_{i=1}^d t^i A_i(b_1-b_2)$,
    by linearity of $F$, we have
    \[ q(t)=\sum_{i=1}^{d} t^i F(A_i(b_1-b_2)),\] so $q(t)$ is a
    polynomial map from $\RR^k$ to $\RR$.

    By our assumptions, the coefficient of $t^r$ in $q(t)$ is not
    zero, hence $q(t)$ is a non-zero polynomial and
    $\lim_{t\to \infty} |q(t)|=\infty$.

    On the other hand, by the definition of $X_j$, we have
    $\rho_t(b_1)-\rho_t(b_2)\in L_j+\bB_1-\bB_1$ for $t\gg 0$.  Since
    $F$ vanishes on $L_j$, we obtain $q(t)\in F(\bB_1-\bB_1)$, which
    is compact set. A contradiction with
    $\lim_{t\to \infty} |q(t)|=\infty$.

    This finishes the proof of the claim.
  \end{proof}

  Let $j=1,\dotsc,n$. By Claim A, for $i=1,\dotsc,d$, we may choose
  $a_{ij}\in L_j^\perp$ with $A_i X_j \subseteq a_{ij}+L_j$.  For
  $x\in X_j$ and $t\in I_\infty$ we have
  \[ \rho_t(x)=\sum_{i=1}^d t^i A_ix\in \sum_{i=1}^d t^i( a_{ij} +L_j)
    = \sum_{i=1}^d t^i a_{ij} + \sum_{i=1}^d t^iL_j. \] Since $L_j$ is
  closed under multiplication by scalars, setting
  $p_j(t)= \sum_{i=1}^d t^i a_{ij}$ we obtain
  \begin{equation}
    \label{eq:2}
    \rho_t(X_j) \subseteq p_j(t)+L_j \text{  and }
    \rho_t(X) \subseteq \bigcup_{j=1}^n (p_j(t)+L_j) \text{ for all } t\in
    \RR.
  \end{equation}
  Notice that each $p_j(t)$ is a proper polynomial map, and takes
  values in $L_j^\perp$.

  \medskip

  Let $\CM$ be the proper polynomial family of multi-cosets
  \[ \CM=\Bigl\{ \bigcup_{j=1}^n (p_j(t)+L_j) : t\in I_\infty
    \Bigr\}. \]

  By (\ref{eq:2}), $\CM$ satisfies clause (1) of the theorem, and by
  the explanation right after (\ref{eq:5}), $\SLm(\CF)=\SLm(\CM)$,
  implying clause (2).  Let us see that clause (3) also holds.

  \medskip

  \noindent\textbf{Claim B.} {\itshape For every $j\!=\!1,\dotsc,n$,
    there exists $c_j\in L_j^\perp$, such that
    $c_j=\lim_{t\to \infty}(a_j(t)-p_j(t))$.}
  \begin{proof}[Proof of Claim B]

    We fix $j=1,\dotsc,n$, and choose $b\in X_j$.  By the definition
    of $X_j$, we have
    \[ \rho_t(b )\in \bB_1+a_j(t)+L_j \text{ for } t\gg 0,\] and by
    (\ref{eq:2}),
    \[ \rho_t(b) \in p_j(t) + L_j.\]

    Since both $a_j(t)$ and $p_j(t)$ take values in $L_j^\perp$, we
    have
    \[
      p_j(t) \in a_j(t)+\pi_j^\perp(\bB_1).
    \]
    Thus $ a_j(t) - p_j(t)\in \pi_j^\perp(\bB_1)$. The set
    $\pi_j^\perp(\bB_1)$ is compact, hence, by o-minimality,
    $\lim_{t\to +\infty} (a_j(t)-p_j(t))$ exists and it belongs to
    $L_j^\perp$, call it $c_j$.  This finishes the proof of Claim
    B.\end{proof}

  Recall that each $a_j(\tau)+L_j$ is a nearest coset to some
  $q_j\in S_{\mathcal Y}(\tau)$. It follows that $p_j(\tau)+ c_j+L_j$
  is also a nearest coset to $q_j$.

  \medskip

  \noindent\textbf{Claim C.} {\itshape For every $j=1\dotsc,n$, we
    have $c_j=0$, namely the coset $p_j(\tau)+L_j$ is a nearest coset
    to $q_j$. }
  \begin{proof}[Proof of Claim C]
    We proceed by reverse induction on $\dim(L_j)$, so we consider
    $j_0\in \{1,\ldots, n\}$ and assume that for all $j$ with
    $\dim(L_j)>\dim(L_{j_0})$, we already know that $c_j=0$, so
    $p_j(\tau)+L_j$ is a nearest coset to $q_j$.

    Since $X=\bigcup_{j=1}^n X_j$, there exists
    $j_1\in \{1,\dotsc,n\}$ such that $q_{j_0}$ lies on
    $\rho_\tau(X_{j_1}^\sharp)$.  From \eqref{eq:2} we conclude
    \begin{equation}
      \label{eq:3}
      q_{j_0}\vdash p_{j_1}(\tau)+L_{j_1}^\sharp,
    \end{equation}
    hence, by the definition of a nearest coset $L_{j_0}\subseteq L_{j_1}$.

    If $L_{j_0} \neq L_{j_1}$, then $\dim(L_{j_1})> \dim(L_{j_0})$, so
    by our induction assumption, the coset
    $p_{j_1}(\tau) + L_{j_1}^\sharp$ is a nearest coset to
    $q_{j_1}$. Thus
    $\mu+p_{j_1}(\tau)+L_{j_1}^\sharp=\mu+a_{j_1}(\tau)+L_{j_1}$,
    hence
    \[ a_{j_0}(t)+L_{j_0} \subseteq \bB_1+ a_{j_1}(t) + L_{j_1},
      \text{ for } t\gg 0. \] It follows that
    \[ \rho_t(X) \subseteq \bigcup_{j=1}^n (a_j(t)+L_j+\bB_1)
      \subseteq \bigcup_{\substack{1\leq j \leq n \\ j\neq j_0}}
      (a_j(t)+L_j+\bB_1),
    \]
    contradicting minimality of $n$.

    Thus, we must have $L_{j_0}= L_{j_1}$.  From the definition of a
    nearest coset and (\ref{eq:3}), we get
    \[ \mu+p_{j_0}(\tau)\ +c_{j_0}+L_{j_0}^\sharp=
      \mu+p_{j_1}(\tau)+L_{j_1}^\sharp.\] Since
    $c_{j_0}\in L_{j_1}^\perp$, and both $p_{j_1}(t)$, $p_{j_0}(t)$
    take values in $L_{j_1}^\perp$, it follows that
    $\lim_{t\to \infty}(p_{j_1}(t) -p_{j_0}(t)) = c_{j_0}$.  Recall
    that both polynomial $p_{j_1}(t)$ and $p_{j_0}(t)$ have zero
    constant terms, so $c_{j_0}=0$.

    This finishes the proof of Claim C.
  \end{proof}

  In order to prove that $\CM$ satisfies clause (3), we need to show
  that for every lattice $\Gamma \subseteq \RR^n$, and
  $\pi:\RR^n\to \RR^n/\Gamma$, the families $\pi(\CF)$ and $\pi(\CM)$
  have the same Hausdorff limits at $\infty$. Equivalently, by
  Proposition~\ref{prop:hlimpi}, it is sufficient to show, for any
  $\tau>\RR$,
\[\st(\rho_\tau(X^\sharp)+\Gamma^\sharp)=\st\Bigl(\bigcup_{j=1}^n
p_j(\tau)+L_j^\sharp+\Gamma^\sharp\Bigr).\] We fix $\tau>\RR$.

By clause (1), we clearly have the left-to-right inclusion. For the
opposite inclusion, for each $j=1,\ldots, n$, by Claim C, there exists
a type $q_j\in S_{\CY}(\tau)$, whose nearest coset is
$p_j(\tau)+L_j$. By Theorem~\ref{thm:main types},
$\st(q_j(\mfR)+\Gamma^\sharp)=\st(p_j(\tau)+L_j+\Gamma^\sharp)$, thus
$\st(\CY+\Gamma^\sharp)\supseteq \st(\bigcup_{j=1}^n
p_j(\tau)+L_j^\sharp+\Gamma^\sharp)$.

This ends the proof of Theorem~\ref{thm:main-dilation}.

\end{proof}

\subsubsection{Hausdorff limits of families of multi-cosets}
\label{sec:hausd-limit-multi}

In this section we describe Hausdorff limits of families of
multi-cosets in real tori. Together with
Theorem~\ref{thm:main-dilation} it provides a complete description of
the Hausdorff limits at $\infty$ of proper families of polynomial
dilations.

\medskip

Let $a_1(t),\ldots, a_n(t)\colon I_\infty\to \mathbb R^m$ be
$\Rom$-definable functions, and let $L_1,\ldots, L_n\subseteq \mathbb R^m$
be linear subspaces.

For $t\in I_\infty$, let $ \CM_t$ be the multi-coset
$\CM_t= \bigcup_{i=1}^n \bigl(a_i(t)+L_i\bigr)$, and let
$\CM=\{\CM_t:t\in I_\infty\}$ be the corresponding family of
multi-cosets.  Notice, we do not assume that $a_i(t)$ are polynomials.

\medskip

Let $\Gamma\subseteq \mathbb R^m$ be a lattice.  We denote by
$\Gamma^n\subseteq (\mathbb R^m)^n$ the lattice obtained by the $n$-fold
cartesian power of $\Gamma$, and as before, for a subspace
$V\subseteq (\mathbb R^m)^n$, we denote by $V^{\Gamma^n}$ the smallest
linear $\Gamma^n$-rational subspace of $(\mathbb R^m)^n$ containing
$V$. For the quotient map $\pi\colon \RR^m\to \RR^m/\Gamma$, we denote
by $\pi^{(n)}$ the quotient map
$\pi^{(n)}\colon (\RR^{m})^n\to (\RR^{m})^n/\Gamma^n$.

\begin{thm}\label{thm:hlom-multi}
  Let $a_1(t),\ldots, a_n(t)$, $L_1,\ldots, L_n$ and $\CM$ be as
  above.  Then, there exists a coset $\bar c+V$ of a linear subspace
  $V\subseteq (\RR^m)^n$, such that for every lattice
  $\Gamma\subseteq \mathbb R^m$, the family of Hausdorff limits of
  $\pi_\Gamma(\CM)$ at $\infty$ is exactly the family
 \[\left\{\pi_\Gamma\Bigl ( \bigcup_{i=1}^n  (d_i+L_i^\Gamma)\Bigr ) : (d_1,\ldots, d_n)\in  \bar c+V^{\Gamma^n}\right\}.\]
\end{thm}

\begin{proof}
  Consider the $\Rom$-definable curve
  $\sigma \colon I_\infty\to (\mathbb R^m)^n$, given by
  $ \sigma(t)=(a_1(t),\ldots, a_n(t))$. Let $p(x)\in S(\RR)$ be the
  unique o-minimal type on $\sigma(t)$ at $\infty$, whose realization
  is the set
  \[ p(\mfR)= \{ \sigma(\tau) \colon \tau \in \mfR, \tau > \RR \}
    \subseteq (\mfR^m)^n.\]

  Let $\bar c+V$ be a nearest coset to $p(x)$.  Since $p$ is a type
  over $\RR$, we have $\bar c=(c_1,\ldots, c_n)$, with each $c_i$ in
  $\mathbb R^m$, and $V\subseteq (\RR^m)^n$ is a subspace.  We claim
  that this coset satisfies the conclusion of the theorem.

  \medskip

  Let $\Gamma \subseteq \RR^m$ be a lattice, and
  $X\subseteq \RR^m/\Gamma$. We let $\pi=\pi_\Gamma$. We denote by
  $\CH_\Gamma$ the set of Hausdorff limits at $\infty$ of $\pi(\CM)$

  Since $\pi$ is $\Gamma$-invariant, it is sufficient to show that
  $X\in \CH_\Gamma$ if and only if
  \[ X=\pi\Bigl ( \bigcup_{i=1}^n (d_i+L_i^\Gamma)\Bigr ) \text{ for
      some } (d_1,\ldots, d_n)\in \bar c+V^{\Gamma^n}+\Gamma^n.\]

  Using Proposition~\ref{prop:hlimpi}, we have that $X\in \CH_\Gamma$
  if and only if
  \[ X=\pi\Bigl(\st \bigcup_{i=1}^n (
    a_i(\tau)+L_i^\sharp+\Gamma^\sharp)\Bigr) \text{ for some }\tau\in
    \mfR \text{ with }\tau >\RR.\] Thus $X\in \CH_\Gamma$ if and only
  if
  \[ X=\pi\Bigl(\bigcup_{i=1}^n \st(
    \alpha_i+L_i^\sharp+\Gamma^\sharp)\Bigr) \text{ for some }
    (\alpha_1,\dots,\alpha_n)\in p(\mfR). \]

  Let $\alpha_1,\dotsc,\alpha_n\in \mfR^n$. By
  Lemma~\ref{lem:fst-cyl}, for every $i=1,\dotsc,n,$ the set
  $\st(\alpha_i+\Gamma^\sharp)$ is non-empty, and for any
  $d_i\in \st(\alpha_i+\Gamma^\sharp)$ we have
  $\st(\alpha_i+L_i^\sharp+\Gamma^\sharp)= d_i + L_i^\Gamma + \Gamma$.

  Also, clearly, for any $d_1,\dotsc,d_n\in \RR^m$ and
  $\alpha_1,\dotsc,\alpha_n\in \mfR^m$ we have
  \[ \bigwedge_{i=1}^n d_i \in \st(\alpha_i +\Gamma^\sharp)\
    \Longleftrightarrow\ (d_1,\dotsc,d_n)\in \st
    \bigl((\alpha_1,\dotsc,\alpha_n)+(\Gamma^n)^\sharp\bigr).\]

  It follows that $X\in \CH_\Gamma$ if and only if
  \[ X=\pi\Bigl(\bigcup_{i=1}^n ( d_i+L_i^\Gamma) \Bigr) \text{ for
      some } (d_1,\dots,d_n)\in
    \st\bigl(p(\mfR)+(\Gamma^n)^\sharp\bigr). \]

  By Theorem~\ref{thm:main types} (over the parameter set
  $\RR$),
  \[\st\bigl(p(\mfR)+(\Gamma^n)^\sharp\bigr) =\st\bigl (\bar c +
  V^\sharp +(\Gamma^n)^\sharp\bigr),\] and by
  Lemma~\ref{lem:st-closed}(2), the set on the right equals
  $\cl(\bar c+V+\Gamma)$.  By Fact~\ref{fact:ratner}, we
  have
\[\st\bigl(p(\mfR)+(\Gamma^n)^\sharp\bigr) =\bar
  c+V^{\Gamma^n}+\Gamma^n.\]

  This finishes the proof of the theorem.

\end{proof}

Theorem \ref{thm:main-dilation} and Theorem \ref{thm:hlom-multi}
immediately yield the definability of the family of Hausodrff limits
at $\infty$ of proper families of polynomial dilations, in the
following sense.
\begin{cor}\label{cor:dil-def}
  Let $\{ \rho_t\colon \RR^k\to \RR^m\colon t\in I_\infty\}$ be a
  proper family of polynomial dilations, and $X\subseteq \RR^k$ an
  $\Rom$-definable set.

  Then there are linear subspaces $L_1,\ldots,L_n\subseteq \RR^n$, and
  there is a coset of a linear space $\bar c+V\subseteq (\RR^m)^n$ such
  that for any lattice $\Gamma\subseteq \RR^m$, the set of Hausdorff
  limits at $\infty$ of the family
  $\{\pi{\co}\rho_t(X) \colon t\in I_\infty\}$ is exactly the family
 \[\left\{\pi\Bigl ( \bigcup_{i=1}^n  (d_i+L_i^\Gamma)\Bigr ) : (d_1,\ldots, d_n)\in  \bar c+V^{\Gamma^n}\right\}.\]

 In particular, it is the projection under $\pi$ of a definable family
 of subsets of $\RR^m$.

\end{cor}

\subsection{Polynomial dilations in unipotent groups}
\label{sec:polyn-dial-unip}
We now consider the case of general unipotent groups.  For the setting
we refer to Section \ref{sec: setting}

\subsubsection{Abelinization of a family of dilations}\label{sec:abel-family-dial}

Following \cite{KSS}*{Section~1.5} we introduce \emph{the
  abelinization} of a family of dilations.

Let $G$ be a unipotent group of dimension $d$, $\Gab=G/[G,G]$ its
abelinization and $\piab\colon G\to \Gab$ the projection map.  Since
$\Gab$ is an abelian unipotent group, we identify it with $(\RR^m,+)$
for $m=\dim(\Gab)$.  We also identify the Lie algebra $\gab$ with
$(\RR^m,+)$, and assume that the exponential map
$\exp\colon \gab\to \Gab$ is the identity map.

Let $d\piab$ be the differential of $\piab$ at the identity $e\in G$.
We have the following commutative diagram with polynomial maps.

\[
  \begin{diagram}
    \node{G} \arrow{e,t}{\piab} \node{\Gab} \\ \node{\mfg}
    \arrow{n,r}{\exp}\arrow{e,b}{d\piab}\node{\mfg=\RR^m}
    \arrow{n,r,<>}{\mathrm{id}}
  \end{diagram}
\]

Let $M_t$ be a polynomial $d\times k$ matrix $M_t$, and
$\{\rho_t\colon \RR_k\to G \,\colon\,t\in I_\infty\}$ the
corresponding polynomial family of dilations.  For $t\in I_\infty$ we
denote by $L_t\colon \RR^k\to \RR^d$ the linear map $x\mapsto
M_tx$. Thus $\rho_t=\exp{\co}L_t$, and the following diagram is
commutative

\[
  \begin{diagram}
    \dgARROWPARTS=10 \dgLABELOFFSET=.5ex \node[2]{G}
    \arrow{e,t}{\piab} \node{\Gab}
    \\
    \node{\RR^k}\arrow{e,b}{L_t}\arrow{ne,t}{\rho_t}
    \arrow{ene,b,..,8}{\rhoab_t} \node{\mfg}
    \arrow{n,r,2}{\exp}\arrow{e,b}{d\piab}\node{\mfg=\RR^m}
    \arrow{n,r,<>}{\mathrm{id}}
  \end{diagram}
\]
with $\rhoab_t=\piab{\co}\rho_t=d\piab{\co}L_t$.

Since $d\piab$ is a linear map, there is a $d\times m$ matrix $D$ such
that $d\piab \colon x \mapsto Dx$. Thus
\[ \rhoab_t\colon x\to D(M_tx)=(DM_t)x.\] It is not hard to see that
$DM_t$ is a polynomial matrix, hence we obtain that the family
$\{\rhoab_t \colon \RR^k\to \Gab\colon t\in I_\infty\}$ is a
polynomial family of dilations, and it is proper if the original
family $\{\rho_t \colon t\in I_\infty\}$ was.

We call the family $\{ \rhoab_t\colon t\in I_\infty\}$ \emph{the
  abelinization of the family} $\{\rho_t \colon t\in I_\infty\}$.

\medskip

We are now ready to prove a strong version of
Theorem~\ref{thm:main-hlin1} for proper polynomial dilations.

\begin{thm}\label{thm:dila-main}
  Let $G$ be a unipotent group,
  $\{ \rho_t\colon \RR^k\to G \colon t\in I_\infty\}$ a proper family
  of polynomial dilations, $X\subseteq \RR^k$ an $\Rom$-definable set,
  and let $\CF$ be the family
  \[ \CF=\{ \rho_t(X) \subseteq G \colon t\in I_\infty\}. \]

  Then, for every lattice $\Gamma\subseteq G$ the following conditions
  are equivalent:
  \begin{enumerate}[(a)]
  \item $L^G=G$ for some $L\in \SLm(\CF)$.
  \item $\pi(\CF)$ converges strongly to $G/\Gamma$ at $\infty$.

  \item $\pi(\CF)$ converges to $G/\Gamma$ at $\infty$.
  \item $G/\Gamma$ is a Hausdorff limit at $\infty$ of the family
    $\pi(\CF)$.
  \end{enumerate}
\end{thm}
\begin{proof}
  By Theorem~\ref{thm:main-hlin1}, $(a)$ snd $(b)$ are equivalent.

  \medskip $(b)\Rightarrow (c)$ and $(c) \Rightarrow (d)$ are obvious.

  \medskip We are left to show that $(d) \Rightarrow (a)$.

  Assume $(d)$ holds. i.e. $G/\Gamma$ is one of the Hausdorff limits
  at $\infty$ of the family $\pi(\CF)$.

  Let $\pi^*$ be the quotient map $\pi^*\colon \Gab\to \Gamab$, and
  $\CFab$ be the family
  \[ \CFab=\{ \rhoab_t(X) \subseteq \Gab \colon t\in I_\infty\}. \]

  Applying abelianization, we obtain that $\Gab/\Gamab$ is a Hausdorff
  limit at $\infty$ of the family $\pi^*(\CFab)$.

  Let $\CM=\{\bigcup_{j=1}^n(p_j(t)+L_j)\colon t\in I_\infty\}$ be a
  proper polynomial family of multi-cosets of $\Gab$, as in
  Theorem~\ref{thm:hlom-multi} applied to the family $\CFab$.  Then
  $\Gab/\Gamab$ is a Hausdorff limit at $\infty$ of $\CM$, and hence
  (e,g. by Lemma~\ref{lem:fst-cyl}(1) and Lemma~\ref{claim:hinfty})
  \[ \Gab= \Gamab+\bigcup_{j=1}^n(a_j+L_j^{\Gamab}), \] for some
  $a_1,\dotsc,a_n\in \Gab$.

  Since $\Gamab$ is a discrete subgroup, it follows that
  $\Gab=L_k^{\Gamab}$ for some $k\in\{1,\dotsc,n\}$.  By
  Fact~\ref{fact:rational-normal}, it follows then that
  $\Gab=L_k^{\Gamma_0}$ for any subgroup $\Gamma_0\subseteq \Gamab$ of
  finite index. Thus $\pi^*(\CM)$ converges strongly to $\Gab/\Gamab$
  at $\infty$, and hence, by the choice of $\CM$, the family
  $\pi^*(\CFab)$ converges strongly at $\infty$ to $\Gab/\Gamab$ as
  well.

  By Corollary~\ref{cor:abel}, $\pi(\CF)$ converges strongly at
  $\infty$ to $G/\Gamma$. Thus $(b)$ holds.

  This finishes the proof of the theorem.
\end{proof}

\begin{rem}
  Theorem~\ref{thm:dila-main} can be compared to
  \cite{KSS}*{Theorem~1.3}. The latter is an equidistribution result
  on measures which are associated to polynomial dilations of real
  analytic curves in nilmanifolds. The set $\SLm$ in our analysis is
  replaced there by (a-priori infinitely many) kernels of
  characters. The equidistribution of the measures implies the
  convergence of the family to $G/\Gamma$, under the appropriate
  assumptions. Our additional input, under the assumption of
  $\Rom$-definability, is the treatment of higher dimensional sets, as
  well as the fact that the sets in $\SLm$ work for all lattices.
\end{rem}

\bibliographystyle{acm}

\end{document}